\documentclass[10pt,a4paper,onefignum,oneeqnum,onetabnum,reqno,final]{siamart251216}
\usepackage[numbers,sort&compress]{natbib}
\usepackage{fancyhdr}

\expandafter\def\expandafter\normalsize\expandafter{%
    \normalsize%
    \setlength\abovedisplayskip{1pt}%
    \setlength\belowdisplayskip{1pt}%
    \setlength\abovedisplayshortskip{1pt}%
    \setlength\belowdisplayshortskip{1pt}%
}

\usepackage{geometry}
\allowdisplaybreaks

\usepackage{array,graphicx,subeqnarray}

\usepackage{caption}
\usepackage{algorithm}
\usepackage{algpseudocode}

\newcommand{\Ex}{ \mathbb{E} }

\def\esssup_#1{\underset{#1}{\mathrm{ess\,sup\, }}}
\def\essinf_#1{\underset{#1}{\mathrm{ess\,inf\, }}}
\def\argmax_#1{\underset{#1}{\mathrm{arg\,max\, }}}
\def\argmin_#1{\underset{#1}{\mathrm{arg\,min\, }}}

\usepackage[english]{babel} 

\usepackage{enumerate}

\newsiamthm{assumption}{Assumption}

\newsiamthm{remark}{Remark}

\usepackage{newpxmath} 

\begin{document}

\title{Continuous-time \lowercase{q}-learning for Markov regime-switching system under Tsallis entropy\thanks{Supported by the National Natural Science Foundation of China (NSFC Grant Nos. 12171086).}}

\date{\vspace{-0.5in}}

\author{Minghui Zhang \thanks{ School of Mathematics, Southeast University, Nanjing, 211189, P. R. China. E-mail: minghuibest@163.com}
\and
Xun Li \thanks{ Department of Applied Mathematics, Hong Kong Polytechnic University, Hong Kong, P. R. China. E-mail: li.xun@polyu.edu.hk}
\and Jie Xiong \thanks{Department of Mathematics and SUSTech International Center for Mathematics, Southern University of Science and Technology, Shenzhen, Guangdong, 518055, China; E-mail: xiongj@sustech.edu.cn}
\and
Xin Zhang \thanks{Corresponding author: School of Mathematics, Southeast University, Nanjing, 211189, P. R. China. E-mail: xzhangseu@seu.edu.cn}
}

\date{}

\maketitle

\smallskip

\noindent \textsc{Abstract}:
This paper studies continuous-time q-learning (the continuous-time counterpart of Q-learning) for a Markov regime-switching system under Tsallis entropy regularization. 
The Tsallis entropy regularization yields an optimal policy distribution that may not necessarily be a Gibbs measure, thereby complicating algorithm design. Furthermore, to address the limited universality of current continuous-time regime-switching reinforcement learning algorithms (often restricted to the exploratory mean-variance framework), this study focuses on continuous-time q-learning for Markov regime-switching systems based on Tsallis entropy, aiming for a more universally applicable continuous-time reinforcement learning method. We establish the martingale characterization of the q-function under Tsallis entropy for continuous-time Markov regime-switching systems. We further design two q-learning algorithms that differ based on whether the Lagrange multiplier can be explicitly derived. We apply these algorithms to the continuous-time exploratory mean-variance portfolio optimization problem in a regime-switching market. Numerical experiments demonstrate the satisfactory performance of our q-learning algorithms.

\vspace{0.1in}

\section{Introduction}
In recent years, the theory and applications of Reinforcement Learning (RL) have made significant advances within continuous-time frameworks. These frameworks have demonstrated particular advantages when addressing complex dynamic systems. As a foundational reinforcement learning algorithm, Q-learning (see \citet{watkins1989learning,watkins1992q}) evaluates the long-term value of state-action pairs using a Q-function and improves policies based on these evaluations. Traditionally, Q-learning has been applied to discrete-time Markov decision processes (MDPs). In this context, the algorithm learns a Q-function that maps state-action pairs to their expected returns, allowing for policy updates by selecting actions that maximize future rewards (\citet{sutton2018reinforcement}). However, applying this traditional approach to continuous-time scenarios presents numerous challenges. The Q-function is a discrete-time concept that reduces to a value function independent of actions in continuous-time settings; thus, it cannot be used to rank or select current actions. \citet{tallec2019making} even argued that “there is no Q-function in continuous-time”. The common approach is to discretize continuous time to obtain a discretized Q-function, then apply existing Q-learning algorithms to address continuous-time problems. Nevertheless, empirical evidence indicates that this method is highly sensitive to the choice of time discretization and performs poorly with small time steps. Instead of discretizing, \citet{kim2021hamilton} adopted an alternative approach, integrating actions as state variables in continuous-time systems. This method ensures that the motion process remains continuously smooth over time with a bounded growth rate. However, they consider only deterministic dynamic systems and discretize the continuous-time problem, thereby requiring actions to be absolutely continuous, which imposes significant limitations in real-world applications.

To address the limitations of traditional Q-learning in continuous-time settings, the pioneering work of \citet{wang2020reinforcement} and \citet{jia2022policyaa,jia2022policymm,jia2023q} established a theoretical foundation for continuous-time reinforcement learning. Their research introduced a series of important theories for situations with continuous state spaces (diffusion processes) and potentially continuous action spaces. \citet{wang2020reinforcement} introduced an exploratory stochastic control framework that incorporates entropy regularization. They achieved this by formulating an entropy-regularized, distribution-valued stochastic control problem for diffusion processes, thereby enabling strategic control rather than predetermined or non-adaptive control.
\citet{jia2022policyaa} focused on policy evaluation (PE), which involves learning the value function of a given random policy and characterizing it as a martingale problem. Building upon this, \citet{jia2022policymm} subsequently addressed the policy gradient (PG) problem, which involves computing the gradient of the learned value function with respect to the current policy. They demonstrated that the PG problem is mathematically reducible to a simpler PE problem, which can then be solved using the martingale approach developed in \citet{jia2022policyaa}. These theoretical contributions have enabled the development of various online and offline Actor-Critic (AC) algorithms for model-free RL in the context of diffusion dynamics. Significantly, many of these derived algorithms provide novel interpretations of, or even recover, well-known algorithms from the discrete-time Markov decision process setting. In a significant advance, \citet{jia2023q} proposed a continuous-time q-learning method. This approach extends the traditional Q-function and Q-learning algorithms to continuous time. The extension is achieved by using a first-order approximation of the advantage function, defined as the temporal difference between the Q-function and the value function. Compared with discrete-time RL algorithms, continuous-time RL methods design policy-iteration rules and loss functions for policy evaluation directly in a continuous-time framework, without requiring full prior discretization. This allows the algorithm to maintain greater stability and robustness during execution, particularly compared to methods that are sensitive to the chosen discretization step size.

The continuous-time RL framework also facilitates the application of advanced mathematical tools and techniques, such as stochastic differential equations (SDEs) and control theory, which are essential for establishing a solid theoretical foundation for the algorithms. In recent years, continuous-time RL theories and algorithms have been extended to address a variety of new problems. For instance, \citet{wang2023reinforcement} proposed an entropy-regularized Actor-Critic RL algorithm to address the optimal execution problem within the continuous-time Almgren-Chriss model. Specifically, the study formulated the optimal execution strategy by constructing a relaxed stochastic control problem under an entropy-regularized mean-variance objective. \citet{dai2023learning} applied RL to Merton's utility maximization problem in an incomplete market, focusing on learning optimal portfolio strategies without requiring knowledge of model parameters. \citet{bo2023optimal} studied the infinite-horizon optimal tracking portfolio problem with capital injection in an incomplete market. Notably, for the case of unknown parameters, they introduced an entropy-regularized exploratory formulation and proposed a continuous-time q-learning algorithm based on a reflected diffusion process to determine the optimal tracking strategy. \citet{han2023choquet} integrated the Shannon regularization term into continuous-time entropy-regularized RL, deriving explicit solutions for optimal strategies within the linear-quadratic (LQ) setting. \citet{WY2024} extended the continuous-time q-learning algorithm to mean-field control problems, in which the integrated q-function and the essential q-function, together with test policies, play crucial roles in their model-free algorithm. \citet{giegrich2024convergence} investigated the global linear convergence of PG methods for exploratory LQ control problems in continuous time. Their work employs geometry-aware gradient descent and proposes a novel algorithm specifically for discrete-time policies.

In the study of continuous-time q-learning, entropy regularization is a significant direction, with Shannon entropy being a commonly used measure. Moreover, \citet{tsallis1988possible} introduced a generalization of the Shannon entropy known as the Tsallis entropy, which offers greater flexibility and robustness in handling learning tasks with diverse policy distributions. Tsallis entropy is particularly advantageous in scenarios involving non-Gaussian, heavy-tailed distributions with compact support. By adjusting the index parameters of Tsallis entropy regularization, it can transform learned optimal policies into various types, providing greater flexibility in handling uncertainty and incentivizing exploration. Notably, \citet{lee2018sparse,lee2019tsallis} extensively studied a class of Markov decision processes (MDPs) with Tsallis entropy. Specifically, \citet{lee2018sparse} introduced causal sparse Tsallis entropy regularization to propose sparse MDPs. Meanwhile, \citet{lee2019tsallis} developed Tsallis MDPs, which use an additional entropy index to provide a unified framework for entropy-regularized RL problems, including the well-known Shannon-Gibbs (SG) entropy. By controlling the entropy index, Tsallis MDPs can generate different categories of optimal policies. Moreover, \citet{donnelly2024exploratory} investigated optimal control problems in models with latent factors, in which the agent controls the distribution of behavior via Tsallis entropy rewards for exploration in both discrete and continuous-time settings. Therefore, this paper considers the Tsallis entropy to encourage exploration.

In many real-world scenarios, phenomena or systems exhibit properties of state transitions or trend changes. Mathematically, these situations are often modeled using Markov chains. In this context, the state transitions of the Markov chain represent changes in market trends. Consequently, Markov chains can effectively simulate random fluctuations in real-world environments, thereby better capturing sudden shifts and uncertainties in dynamic settings. In light of this, \citet{zy2003markowitz} investigated the continuous-time mean-variance (MV) problem with regime-switching, but their work did not integrate it with reinforcement learning. More recently, \citet{wang2020reinforcement} proposed the exploratory mean-variance (EMV) framework. This approach replaces deterministic strategies with probability distributions and incorporates entropy regularization into the portfolio optimization problem. They solved for the optimal investment strategy, finding it to be a probability distribution over the control space rather than a deterministic control function. Recent research by \citet{WL2024} examined the application of reinforcement learning to the continuous-time mean-variance portfolio problem in a regime-switching market. They formulated the problem as a partially observable Markov decision process (POMDP) with two unobservable market regimes and addressed the partial observability using the Wonham filter. However, their approach relies on several key assumptions: state-independent volatility and that market parameters are either known or estimated from data. Specifically, the martingale loss function employed by \citet{WL2024} requires knowledge of the true market parameters, which is impractical in real-world applications. To address these limitations, \citet{chen2025exploratory} investigated the exploratory mean-variance with regime-switching (EMVRS) problem. Their model extends the existing framework by allowing for state-dependent volatility. Crucially, their parameter update scheme uses the orthogonality-condition loss derived from the optimal EMVRS value function. The advantage of this loss function is that it does not depend on the true market parameters, allowing the RL algorithm to learn and infer the latent market parameters directly from the training data. Simulation studies further confirm that, through this method, the parameters of the proposed RL model converge to their corresponding true values, even when initialized from randomly selected starting points.

Continuous-time q-learning is underexplored, particularly regarding entropy regularization. Current RL algorithms predominantly employ Shannon entropy because the resulting regularized and normalized optimal policy is essentially Gaussian. The ability to explicitly determine its mean and variance greatly simplifies subsequent algorithm design. However, a major challenge with general entropy regularization methods is that the regularized and normalized optimal policy may not have a determined distribution form, which fundamentally constrains the design of policy-based algorithms. To sidestep this issue, some research has shifted toward value-based algorithms, such as temporal-difference (TD) learning and orthogonality-condition (OC) learning (martingale-loss learning), as explored in \citet{chen2025exploratory}. Nonetheless, policy-based algorithms must directly determine the optimal policy. Furthermore, studies on continuous-time q-learning incorporating regime-switching mechanisms remain limited. Therefore, this paper introduces Tsallis entropy to encourage exploration under a regime-switching framework and develops a corresponding continuous-time q-learning method. Our framework extends the work presented in \citet{WY2024} and \citet{chen2025exploratory}. Notably, our model reduces to that of \citet{chen2025exploratory} when the Tsallis entropy order $p$ is set to 1. Furthermore, it is worth emphasizing that our study aims to develop a more general continuous-time reinforcement learning method.

The remainder of this paper is organized as follows. In Section 2, we formulate an entropy-regularized, exploratory RL framework for continuous-time, continuous-space Markov regime-switching systems under the Tsallis entropy, and present some useful preliminary results. In Section 3, we derive the corresponding q-function and establish its fundamental martingale properties, showing that the optimal policy is characterized by the q-function's dependence on a Lagrange multiplier. In Section 4, we design the corresponding q-learning algorithms for two distinct cases, depending on whether the Lagrange multiplier is known or unknown. In Section 5, we present a practical example of the EMVRS problem, including the parametric forms of the optimal value and the q-function, and demonstrate satisfactory convergence in numerical experiments. Finally, Section 6 presents our conclusion. The Appendix contains supplementary materials and detailed proofs for the statements presented in the main text.

\section{Problem Formulation and Preliminaries}\label{sec:Markovmodel}

\subsection{Exploratory formulation of the reinforcement learning problem}
For a fixed time horizon $T>0$, let $(\Omega, \mathcal{F}, \mathbb{F}, \mathbb{P})$ be a filtered probability space with the filtration $\mathbb{F}:=\{\mathcal{F}_s\}_{s \in[0, T]}$ generated by a standard $n$-dimensional Brownian motion $W=\{W_s\}_{s \in[0, T]}$ and a continuous-time irreducible Markov chain $\alpha=\{\alpha_s\}_{s \in[0, T]}$ with a finite state space $\mathcal{M}=\{1, 2, \ldots, L\}$ such that $W_s$ and $\alpha_s$ are independent of each other. We assume the generator of the Markov chain $\alpha$ is $Q = (q_{ij})_{L \times L}$. Let $\Phi_j(s)$ be the number of jumps into state $j$ up to time $s$ and set $\widetilde{\Phi}_j(s)\triangleq \Phi_j(s)-\int_0^s \sum_{i\neq j}^L q_{ij} 1_{\{\alpha_{u-=i}\}}(u) \mathrm{d}u$. Then for each $j\in \mathcal{M},$ the process $\widetilde{\Phi}_j(s)$ is an $(\mathbb{F}, \mathbb{P})$-martingale. For any $i, j \in \mathcal{M}$ and $s^\prime \geq s \geq 0$, we further let $p_{i j}(s^\prime, s)$ denote the transition probability $P(\alpha_{s^\prime}=j | \alpha_s=i)$.

We consider the following controlled Markov regime-switching \textit{state} process:
\begin{align}\label{eq:abstractmodel}
\mathrm{d} X_s^{a}=b(s, X_s^{a}, \alpha_s, a_s) \mathrm{d} s+\sigma(s, X_s^{a}, \alpha_s, a_s) \mathrm{d} W_s,\,\, s \in[t, T] ; \,\, X_t^a=x,\alpha_t=i,
\end{align}
where $a_s$ represents the agent's action or control at time $s$ and $a=\{a_s\}_{s\in [t, T]}$ is an $\mathbb{F}$-adapted process taking values on the action space $\mathcal{A}$. Here, $b:[0, T] \times \mathbb{R}^d \times \mathcal{M}\times \mathcal{A} \mapsto \mathbb{R}^d$, $\sigma:[0, T] \times \mathbb{R}^d \times \mathcal{M} \times \mathcal{A} \mapsto \mathbb{R}^{d \times n}$ are assumed to be measurable functions. 
For each initial time-state-regime triple $(t, x, i) \in [0, T) \times \mathbb{R}^d \times \mathcal{M}$, the agent is interested in the stochastic control problem of finding the optimal control $a=\{a_s\}_{s\in [t,T]}$ that maximizes the following expected discounted total reward: 
\begin{align}\label{eq:objective}
\Ex\left[\int_t^T e^{-\beta(s-t)} r(s, X_s^{a}, \alpha_s, a_s) \mathrm{d} s+ e^{-\beta(T-t)}h(X_T^{a}) | X_t^{a}=x, \alpha_t=i\right],
\end{align}
where $r$ is the running reward function, $h$ is the terminal reward function applied at the end of the period $T$, and $\beta  \geq 0$ is a constant discount factor.

To ensure that the above stochastic control problem 
is well-posed, we impose the following assumptions:
\begin{assumption}\label{assumption}
The following conditions for the state dynamics and reward functions hold true:
\begin{itemize}
\item[(i)] $b, \sigma, r, h$ are all continuous functions in their respective arguments;
\item[(ii)] $b, \sigma$ are uniformly Lipschitz continuous in $x$, i.e., for $\varphi \in\{b, \sigma\}$, there exists a constant $C>0$ such that
\begin{align*}
|\varphi(t, x, i, a)-\varphi\left(t, x^{\prime}, i, a\right)| \leq C|x-x^{\prime}|, \;\; \forall(t, i, a) \in[0, T] \times \mathcal{M} \times \mathcal{A}, \;\; \forall x, x^{\prime} \in \mathbb{R}^d ;
\end{align*}
\item[(iii)] $b, \sigma$ have linear growth in $x$, i.e., for $\varphi \in\{b, \sigma\}$, there exists a constant $C>0$ such that
\begin{align*}
|\varphi(t, x, i, a)| \leq C(1+|x|), \quad \forall(t, x, i, a) \in[0, T] \times \mathcal{M} \times \mathbb{R}^d \times \mathcal{A};
\end{align*}
\item[(iv)] $r$ and $h$ have polynomial growth in $(x, i, a)$ and $x$ respectively, i.e., there exist constants $C>0$ and $c\geq 1$ such that
\begin{align*}
|r(t, x, i, a)| \leq C(1+|x|^{c}+|a|^{c}), \; |h(x)| \leq C(1+|x|^{c}),\; \forall(t, x, i, a) \in[0, T] \times \mathcal{M} \times \mathbb{R}^d \times \mathcal{A}.
\end{align*}
\end{itemize}
\end{assumption}

If the coefficients $b, \sigma, r, h$ in \eqref{eq:abstractmodel}-\eqref{eq:objective} are fully known, classical methods such as the dynamic programming principle (DPP) or stochastic maximum principle can be used to solve the above optimal control problem. However, in practice, the decision maker (agent) often has limited or no knowledge of the environment (i.e., $b,\sigma,r,h$ are unknown). RL offers an efficient method to learn optimal control for \eqref{eq:abstractmodel}-\eqref{eq:objective} in such unknown environments by repeatedly taking actions and interacting with the environment (trial-and-error procedure).

We now present the RL formulation of the problem to be studied in this paper. In what follows, we write 
$\mathcal{P}(\mathcal{A})$ for the set of all probability densities defined on $\mathcal{A}$, i.e., $\mathcal{P}(\mathcal{A}) = \{ \pi \in L^1(\mathcal{A}) \mid \pi(a) \geq 0, \text{ a.e. on } \mathcal{A},\,\, \int_{\mathcal{A}} \pi(a)\mathrm{d}a = 1 \}$. Similar to the approach in \citet{jia2022policyaa,jia2022policymm,jia2023q}, we model exploration by sampling the action or control $a$ at each trial from a probability density  $\pi = \{\pi(a|t, X_t, \alpha_t)\}_{t \geq 0}$ defined over $\mathcal{A}$. The sampled probability density $\pi$ is often referred to as the agent's stochastic policy.  
Fix an initial time-state-regime triple $(t, x, i)$, consider the following SDE
\begin{align}\label{eq:dX_s_pi}
\mathrm{d} X_s^{\pi}=b(s, X_s^{\pi}, \alpha_s, a_s^{\pi}) \mathrm{d} s+\sigma(s, X_s^{\pi}, \alpha_s, a_s^{\pi}) \mathrm{d} W_s,\,\, s \in[t, T] ; \,\, X_t^\pi=x,\alpha_t=i,
\end{align}
where the solution to \eqref{eq:dX_s_pi}, $X^{\pi}=\{X_s^{\pi}, t  \leq s  \leq T\}$, is the sample state process corresponding to $a^{\pi}$. To encourage exploration, we introduce the Tsallis entropy with order $p\geq 1$ as a regularizer, and the objective functional is defined by
\begin{align}\label{eq:value-fun(X_s_pi)}
J^0(t, x, i; \pi)= & \Ex\bigg[\int_t^T e^{-\beta(s-t)}\big[r(s, X_s^{\pi},\alpha_s,a_s^{\pi})+\gamma l_p(\pi(a\mid s,X_s, \alpha_s)))\big] \mathrm{d} s +e^{-\beta(T-t)} h(X_T^{\pi}) \mid X_t^{\pi}=x,\alpha_t^{\pi}=i\bigg],
\end{align}
and $\gamma\geq 0$ stands for \textit{the temperature parameter}. Here, the Tsallis entropy with order $p \geq 1$ is defined by
\begin{align}\label{eq:Tsallis-entropy}
l_p(z)=
\left\{\begin{array}{ll}
\dfrac{1}{p-1}(1-z^{p-1}), & p>1, \\
-\ln z, & p=1,
\end{array}\right. \quad z \in \mathbb{R}_{+}.
\end{align}
However, the representation \eqref{eq:dX_s_pi}-\eqref{eq:value-fun(X_s_pi)} cannot be applied to derive the exploratory HJB equation directly from the point of view of DPP. For this purpose, it is necessary to present a relaxed formulation of the control problem.
Thus, \textit{the exploratory state dynamics} is given by 
\begin{align}\label{eq:state-Exp}
\mathrm{d} \widetilde{X}_s^\pi=\widetilde{b}(s, \widetilde{X}_s^\pi,\alpha_s,\pi(\cdot| s, \widetilde{X}_s^\pi, \alpha_s)) \mathrm{d} s+\widetilde{\sigma}(s, \widetilde{X}_s^\pi,\alpha_s,\pi(\cdot| s, \widetilde{X}_s^\pi, \alpha_s)) \mathrm{d} W_s,
\end{align}
where the coefficients $\widetilde{b}(\cdot, \cdot, \cdot, \cdot)$ and $\widetilde{\sigma}(\cdot, \cdot, \cdot, \cdot)$ are defined by
$$
\widetilde{b}(s, x, i, \pi):=\int_{\mathcal{A}} b(s, x, i, a) \pi(a| s, x, i) \mathrm{d} a, \,\, \widetilde{\sigma}(s, x, i, \pi):=\sqrt{\int_{\mathcal{A}}\sigma \sigma^{\top}(s, x, i, a) \pi(a| s, x, i)\mathrm{d} a}.
$$
Rigorously, by extending the Markovian projection theory (e.g., \citet{brunick2013mimicking}, Theorem 3.6) to regime-switching jump-diffusions (see, e.g., \citet{bentata2012mimicking}), or equivalently by applying the projection conditionally on the sample paths of the independent Markov chain $\alpha$, it is guaranteed that the joint state-regime pair $(X_s^\pi, \alpha_s)$ and its relaxed counterpart $(\tilde{X}_s^\pi, \alpha_s)$ share the exact same marginal probability distribution for each $s \in [t, T]$.
Consider the following objective value function for the policy $\pi$ 
\begin{align}\label{eq:value-fun(MC)}
J(t, x, i; \pi)=\Ex\bigg[\int_t^T e^{-\beta(s-t)} \tilde{\mathcal{R}}(s, \tilde{X}_s^{\pi}, \alpha_s, \pi(a | s, \tilde{X}_s^{\pi}, \alpha_s)) \mathrm{d} s+e^{-\beta(T-t)} h(\tilde{X}_T^{\pi}) | \tilde{X}_t^{\pi}=x,\alpha_t=i\bigg],
\end{align}
where $\tilde{\mathcal{R}}(s, x, i, \pi)=\int_{\mathcal{A}}[r(s, x, i, a)+\gamma l_p(\pi(a| s, x, i))] \pi(a| s, x, i) \mathrm{d} a.$
The task of RL is to find the optimal policy $\pi^*$ and its value function $J(t,x,i;\pi^*)$ such that
\begin{align}\label{eq:opt-value-fun}
V(t, x, i):=J(t,x,i;\pi^*)=\sup_{\pi \in \Pi} J(t, x, i; \pi),
\end{align}
where $\Pi$ stands for the set of admissible (stochastic) policies on $\mathcal{A}$ and $V(t, x, i)$ is called the optimal value function. The following provides the precise definition of admissible (feedback) policies. 




\begin{definition}\label{def:admissible}
A policy $\pi=\pi(\cdot | \cdot, \cdot, \cdot)$ is called admissible if
\begin{itemize}
\item[(i)] $\pi(\cdot | t, x, i) \in \mathcal{P}(\mathcal{A}), supp \,\pi(\cdot | t, x, i)=\mathcal{A}$ for every $(t, x, i) \in[0, T] \times \mathbb{R}^d \times \mathcal{M}$ and $\pi(a | t, x, i):(t, x, i, a) \in [0, T] \times \mathbb{R}^d \times \mathcal{M}\times \mathcal{A} \mapsto \mathbb{R}$ is measurable;
\item[(ii)] the SDE \eqref{eq:state-Exp} has a unique solution for initial $(t, x, i)\in [0,T]\times\mathbb{R}^d \times \mathcal{M}$;
\item[(iii)] $\boldsymbol{\pi}(a | t, x, i)$ is continuous in $(t, x)$ and uniformly Lipschitz continuous in $x$ in the total variation distance, i.e., for any $i\in \mathcal{M}, \lim_{(t^{\prime}, x^{\prime}) \rightarrow (t, x)}\int_{\mathcal{A}}|\pi(a | t, x, i)-\pi(a | t^{\prime}, x^{\prime}, i)| \mathrm{d} a =0$, and there is a constant $C>0$ independent of $(t, x, i)$ such that
\begin{align*}
\int_{\mathcal{A}}|\pi(a | t, x, i)-\pi(a | t, x^{\prime}, i)| \mathrm{d} a \leq C|x-x^{\prime}|, \quad \forall x, x^{\prime} \in \mathbb{R}^d, \quad i\in \mathcal{M};
\end{align*}
\item[(iv)] For any given $n>0$, the entropy of $\pi$ and its $n$-th moment have polynomial growth in $x$, i.e., for any $(t, x, i) \in[0, T] \times\mathbb{R}^d \times \mathcal{M}$, there are constants $C=C(n)>0$ and $c=c(n) \geq 1$ such that
\begin{align*}
    \int_{\mathcal{A}}|a|^{n}  \pi(a | t, x, i) \mathrm{d} a<C(1+|x|^{c}), \quad   \int_{\mathcal{A}} l(\pi(a | t, x, i)) \pi(a | t, x, i) \mathrm{d} a<C(1+|x|^{c}).
\end{align*}
\end{itemize}
\end{definition}

Under the above conditions, we obtain the following result. Since the proof is similar to that of \citet{jia2022policymm}, it is omitted here.

\begin{lemma}\label{lem:growth-condition}
Let Assumptions \ref{assumption} hold and $\boldsymbol{\pi}$ be a given admissible policy. Then the $S D E$ \eqref{eq:state-Exp} admits a unique strong solution that   
satisfies the following growth condition
\begin{align*}
    \Ex\bigg[\max _{t \leq s \leq T}|\tilde{X}_s^\pi|^n | \tilde{X}_t^\pi=x, \alpha_t=i\bigg] \leq C\left(1+|x|^n\right), \quad n\geq 2,
\end{align*}for some constant $C=C(n,i)$. Finally, the expected payoff \eqref{eq:value-fun(MC)} is finite.
\end{lemma}

\subsection{The exploratory HJB equation}\label{subsection:HJB}
By the DPP, the optimal value function defined in \eqref{eq:opt-value-fun} satisfies the following exploratory HJB equation:
\begin{equation}\label{eq:HJB-fun(Exp)}
\left\{\begin{array}{l}\begin{aligned}
 &\sup_{\pi \in \mathcal{P}(A)} \int_{\mathcal{A}}\bigg[H(t, x, i, a, V)+\sum_{j=1}^L q_{ij}V(t, x, j)+\gamma l_p(\pi(a| t, x, i)\bigg]\pi(a\mid t, x, i)\mathrm{d}a =0,\\
&V(T, x, i)=h(x),
\end{aligned}
\end{array}
\right.
\end{equation}
where
\begin{align}\label{eq:Hamiltonian}
H(t, x, i, a, V)=V_t(t,x,i)+b(t, x, i, a)\circ V_{x}(t, x, i)+\frac{1}{2}\sigma \sigma^\top(t, x, i, a)\circ V_{x x}(t, x, i) + r(t, x, i, a)-\beta V(t,x,i).
\end{align}
Here “$\circ$” denotes the inner product.

Now, we first present a lemma that is crucial for finding the optimal feedback policy and establishing the policy improvement theorem. 
    \begin{lemma} \label{lem:unique-maximizer}
        Let $\gamma>0$, $p \geq 1$ and suppose that the function $g(t,x,i,a): [0, T] \times \mathbb{R}^d \times \mathcal{M}\times \mathcal{A} \mapsto \mathbb{R}$ satisfies $0<Z:=\int_{\mathcal{A}}\exp\bigg\{\dfrac{1}{\gamma}g(t,x,i,a)\bigg\}\mathrm{d}a<\infty$. Then, the following optimization problem
        \begin{align}
            \sup_{\pi\in \mathcal{P}(\mathcal{A})} g_p(\pi):=\int_{\mathcal{A}}\left[g(t,x,i,a)+\gamma l_p(\pi(a))\right]\pi(a)\mathrm{d}a
        \end{align}
        admits a unique maximizer
        \begin{align}\label{eq:unique-maximizer-3}
           \pi_p^*(a|t,x,i)=\left\{\begin{array}{cl}
            \left(\frac{p-1}{p \gamma}\right)^{\frac{1}{p-1}}\left(g(t,x,i,a)+\frac{\gamma}{p-1}+\psi_p(t,x,i)\right)_{+}^{\frac{1}{p-1}}, & p>1, \\
             \frac{\exp \left(\frac{1}{\gamma}g(t,x,i,a)\right)}{\int_{\mathcal{A}} \exp \left(\frac{1}{\gamma}g(t,x,i,a')\right) \mathrm{d} a'}, & p=1.
        \end{array}\right.
        \end{align}
  Here,  $\psi_p(t,x,i)$ is called a normalizing function and is determined by
        \begin{align}\label{eq:unique-maximizer-1}
            \int_{\mathcal{A}}\left(\frac{p-1}{p \gamma}\right)^{\frac{1}{p-1}}\left(g(t,x,i,a)+\frac{\gamma}{p-1}+\psi_p(t,x,i)\right)_{+}^{\frac{1}{p-1}} \mathrm{d} a=1.
        \end{align} Furthermore, we also have 
        $\pi_1^*(a|t,x,i)=\lim_{p\downarrow 1}\pi_p^*(a|t,x,i)$. 
    \end{lemma}

\begin{proof}

We introduce a scalar Lagrange multiplier $\psi_p(t, x, i)$ 
to enforce the constraint $\int_{\mathcal{A}} \pi(a)\mathrm{d}a = 1$ and Karush-Kuhn-Tucker (KKT) multiplier $\xi_p(t, x, i, a)$ to enforce the constraint $\pi(a)\geq 0$. The Lagrangian is formulated as
        \begin{equation*}
            \begin{aligned}
                \mathcal{L}(\pi, \psi_p, \xi_p):=& g_p(\pi) + \psi_p(t,x,i) \big( \int_{\mathcal{A}} \pi(a) \mathrm{d}a - 1 \big) + \int_{\mathcal{A}} \xi_p(t,x,i,a) \pi(a) \mathrm{d}a\\
                =& \int_{\mathcal{A}} \bigg[ \big( g(t,x,i,a) + \gamma l_p(\pi(a)) + \psi_p(t,x,i) + \xi_p(t,x,i,a) \big) \pi(a) \bigg] \mathrm{d}a - \psi_p(t,x,i).
            \end{aligned}
        \end{equation*}
In the following, we begin by considering $p>1$.  From the theory of the infinite-dimensional extremum problem, the first G{\^a}teaux derivative of $\mathcal{L}(\pi, \psi_p, \xi_p)$ with respect to $\pi$ must vanish at the optimal policy $\pi_p^*(\cdot|t,x,i)$ along any admissible variation $h \in L^1(\mathcal{A})$. Thus, we have
$$
 g(t,x,i,a) + \frac{\gamma}{p-1} + \psi_p(t,x,i) + \xi_p(t,x,i,a) - \frac{\gamma p}{p-1} \pi_p^*(a|t,x,i)^{p-1} = 0, \quad \text{a.e. } a \in \mathcal{A}.
$$
Solving algebraically for $\pi_p^*(a|t,x,i)$ yields the parametric expression:
\begin{equation}\label{eq:parametric expression}
    \pi_p^*(a|t,x,i) = \left( \frac{p-1}{p \gamma} \right)^{\frac{1}{p-1}} \left( g(t,x,i,a) + \frac{\gamma}{p-1} + \psi_p(t,x,i) + \xi_p(t,x,i,a) \right)^{\frac{1}{p-1}}.
\end{equation}
As $\xi_p(t,x,i,a)$ is the KKT multiplier for the non-negativity constraint, it satisfies $\xi_p(t, x, i, a) \geq 0$ and the complementary slackness condition $\xi_p(t,x,i,a) \pi_p^*(a|t,x,i) = 0$.
Thus, we have the following two cases:
\begin{enumerate}
\item [(i)] If 
            $\pi_p^*(a|t,x,i)>0$, 
            the complementary slackness requires $\xi_p(t,x,i,a) = 0$. Consequently, we have $g(t,x,i,a) + \frac{\gamma}{p-1} + \psi_p(t,x,i) > 0.$ Hence, $\xi_p(t,x,i,a)=0=\big( -(g(t,x,i,a) + \frac{\gamma}{p-1} + \psi_p(t,x,i)) \big)_+$;
\item [(ii)] If 
            $\pi_p^*(a|t,x,i) = 0$, i.e., $g(t,x,i,a) + \frac{\gamma}{p-1} + \psi_p(t,x,i) + \xi_p(t,x,i,a) = 0$. Thus, the KKT multiplier $\xi_p(t,x,i,a) \geq 0$ implies $-\big(g(t,x,i, a) + \frac{\gamma}{p-1} + \psi_p(t,x,i)\big)=\xi_p(t,x,i,a) \geq 0$. Therefore, we also have $\xi_p(t,x,i,a)=\big( -(g(t,x,i,a) + \frac{\gamma}{p-1} + \psi_p(t,x,i)) \big)_+$.
\end{enumerate}
Thus, substituting $\xi_p(t,x,i,a) = \big( -(g(t,x,i,a) + \frac{\gamma}{p-1} + \psi_p(t,x,i)) \big)_+$ back into \eqref{eq:parametric expression}, we obtain
        $$
        \pi_p^*(a|t,x,i) = \bigg( \frac{p-1}{p \gamma} \bigg)^{\frac{1}{p-1}} \bigg( g(t,x,i,a) + \frac{\gamma}{p-1} + \psi_p(t,x,i) \bigg)_{+}^{\frac{1}{p-1}}
        , \quad \text{a.e. } a \in \mathcal{A}.
        $$
The scalar Lagrange multiplier $\psi_p(t,x,i)$ is chosen to ensure that $\pi_p^*(a|t,x,i)$ is a probability density, i.e., it is determined by \cref{eq:unique-maximizer-1}. By the same method as above, we can also derive the expression of $\pi_1^*(a|t,x,i)$ as shown in \eqref{eq:unique-maximizer-3}.

Next, we shall prove the fact that  $\pi_1^*(a|t,x,i)=\lim_{p\downarrow1}\pi_p^*(a|t,x,i)$. 
It is easy to verify that
\begin{align}\label{eq:log-pi-p}
    \log \pi_p^*(a|t,x,i)
    &=\frac{\log\left(1+\frac{p-1}{\gamma}\big(g(t,x,i,a)+\psi_p(t,x,i)\big)\right)_+}
      {p-1}-\frac{\log p}{p-1}.
\end{align}
Thus, if we can prove that $\lim_{p\downarrow 1}\psi_p(t,x,i)=\gamma(1-\log Z)=\gamma(1-\log \int_{\mathcal{A}}\exp\{\dfrac{1}{\gamma}g(t,x,i,a))\} \mathrm{d} a$, then by observing that
$\lim_{y\downarrow0}\frac{\log(1+y)}{y}=1, \lim_{p\downarrow 1}\frac{\log p}{p-1}=1$, we can show that $\lim_{p\downarrow 1}\log\pi_p^*(a|t,x,i)=\dfrac{1}{\gamma}g(t,x,i,a)-\log Z$.
 Therefore, we have $\lim_{p\downarrow 1}\pi_p^*(a|t,x,i)=\dfrac{\exp \left(\frac{1}{\gamma}g(t,x,i,a)\right)}{\int_{\mathcal{A}} \exp \left(\frac{1}{\gamma}g(t,x,i,a')\right) \mathrm{d} a'}=\pi_1^*(a|t,x,i)$.

In the following, we prove that $\lim_{p\downarrow 1}\psi_p(t,x,i)=\gamma(1-\log Z)$. We first show that $\psi_p(t,x,i)$ is bounded as $p\downarrow 1$. Observing that
\begin{align}\label{eq:pip*bound}
    \pi_p^*(a|t,x,i)&=\bigg( \frac{1}{p} \bigg[1+\frac{p-1}{\gamma}\big(g(t,x,i,a)+\psi_p(t,x,i)\big)\bigg]_+\bigg)^{\frac{1}{p-1}}\leq \big[\exp\{\frac{p-1}{\gamma}\big(g(t,x,i,a)+\psi_p(t,x,i)\big)\}\big]^{\frac{1}{p-1}} \notag\\
    &=\exp\{\frac{1}{\gamma}\big(g(t,x,i,a)+\psi_p(t,x,i)\big)\},
\end{align}
where we use the fact that $p>1$ and $(1+y)_{+}\leq e^y$.
Integrating the above inequality on both sides over $\mathcal{A}$ yields $1\leq e^{\psi_p(t,x,i)/\gamma}Z$. Therefore, $\psi_p(t,x,i)\geq -\gamma\log Z$.

On the other bound, we let $\mathcal{A}_N=\{a\in\mathcal{A}:|g(t,x,i,a)|\leq N\}\subset\mathbb{R}^d$ and $m_N$ be the Lebesgue measure of $\mathcal{A}_N$. Since $g(t, x, i, a) \geq-N$ on $\mathcal{A}_N$, we have
\begin{align*}
    1=\int_{\mathcal{A}} \pi_p^*(a|t,x,i)\mathrm{d}a
    \geq
    \int_{\mathcal{A}_N}\pi_p^*(a|t,x,i)\mathrm{d}a
    \geq
    m_N\bigg( \frac{1}{p} \big(1+\frac{p-1}{\gamma}\big(-N+\psi_p(t,x,i)\big)\big)_+\bigg)^{\frac{1}{p-1}}.
\end{align*}
Consequently, we obtain $\psi_p(t,x,i)\leq \gamma \dfrac{pm_N^{1-p}-1}{p-1}+N\to \gamma(1-\log m_N)+N$ as $p\downarrow 1$. Thus, by combining the obtained lower bound, we prove that $\psi_p(t,x,i)$ is bounded as $p\downarrow 1$. Furthermore, from \cref{eq:pip*bound}, we can also prove that $\pi_p^*(a|t,x,i)\leq C\exp\{\frac{1}{\gamma}g(t,x,i,a)\}$ for some constant $C>0$ and all $p$ close to $1$.

Let $p_n$ be an arbitrary sequence such that $p_n \downarrow 1$. By boundedness of $\psi_{p_n}(t,x,i)$, we can select a subsequence (still denoted as $\psi_{p_n}(t,x,i)$) and some $\bar{\psi}(t,x,i) \in \mathbb{R}$, such that $\psi_{p_n}(t,x,i) \longrightarrow \bar{\psi}(t,x,i)$.
Thus, from \eqref{eq:log-pi-p}, we have
\begin{align}\label{eq:log-pi-p-limit}
    \lim_{n \to \infty} \pi_{p_n}^*(a|t,x,i) = \exp\left\{\frac{1}{\gamma}\big(g(t,x,i,a)+\bar{\psi}(t,x,i)\big) - 1\right\}.
\end{align}
Therefore, by the dominated convergence theorem, we have
\begin{align*}
    1=\lim_{n \to \infty} \int_{\mathcal{A}} \pi_{p_n}^*(a|t,x,i) \mathrm{d} a = \int_{\mathcal{A}} \lim_{n \to \infty} \pi_{p_n}^*(a|t,x,i) \mathrm{d} a = \int_{\mathcal{A}} \exp\left\{\frac{1}{\gamma}\big(g(t,x,i,a)+\bar{\psi}(t,x,i)\big) - 1\right\} \mathrm{d} a ,
\end{align*}
which further implies that $\bar{\psi}(t,x,i)=\gamma(1-\log Z)$. 
Hence, we can conclude that as $p \downarrow 1$, the limit of $\psi_p(t,x,i)$ exists and equals $\gamma(1-\log Z)$.

The uniqueness of the maximizer follows directly from the fact that $\mathcal{P}(\mathcal{A})$ is a convex set and $g_p(\pi)$ is a strictly concave functional on $\mathcal{P}(\mathcal{A})$, and $p=1$ is interpreted by the limit as $p \to 1$.
    \end{proof}

We can now apply \cref{lem:unique-maximizer} to determine the maximizer of the exploratory HJB equation \cref{eq:HJB-fun(Exp)}. Let $g(t,x,i,a)=H(t, x, i, a, V)+\sum_{j=1}^L q_{i j} V(t, x, j)$, then, for $p\geq 1$, we can obtain the maximizer of the exploratory HJB equation \cref{eq:HJB-fun(Exp)} as:
\begin{align}\label{eq:HJB-fun(Exp)-maximizer}
    \pi_p^*(a|t,x,i)=\left(\frac{p-1}{p \gamma}\right)^{\frac{1}{p-1}}\bigg(H(t, x, i, a, V)+\sum_{j=1}^L q_{i j} V(t, x, j)+\frac{\gamma}{p-1}+\psi_p(t,x,i)\bigg)_{+}^{\frac{1}{p-1}},
\end{align}
where 
$\psi_p(t,x,i)$ is determined by $\int_{\mathcal{A}}\pi_p^*(a|t,x,i)\mathrm{d}a=1$. 

Next, we will establish the policy improvement based on the maximizer \eqref{eq:HJB-fun(Exp)-maximizer} of the exploratory HJB equation.   
Let us first recall that the objective function $J(t, x, i;\pi)$ for a fixed admissible policy $\pi$ is defined by \eqref{eq:value-fun(MC)}. Thus, if the objective function $J(\cdot, \cdot, i ; \pi) \in C^{1,2}([0, T) \times \mathbb{R}^d) \cap C([0, T] \times \mathbb{R}^d)$ for all $i\in \mathcal{M}$, it satisfies the following PDE:
\begin{align}\label{eq:obj-fun(general)}
 \int_{\mathcal{A}}\bigg[H(t, x, i, a, J(t,x,i;\pi))+\sum_{j=1}^L q_{ij}J(t, x, j;\pi)+\gamma l_p(\pi(a\mid t, x, i))\bigg] \pi(a\mid t, x, i) \mathrm{d} a=0,
\end{align}
with the terminal condition $J(T, x, i;\pi)=h(x)$.

\begin{theorem}[Policy Improvement Iteration]\label{thm:policy-improve}
For any $\pi \in \Pi$, assume that the objective function $J(\cdot, \cdot, i ; \pi) \in C^{1,2}([0, T) \times \mathbb{R}^d) \cap C([0, T] \times \mathbb{R}^d)$ 
for all $i\in \mathcal{M}$, and for $p>1$, there exists a function $\psi(t, x, i;\pi):[0, T] \times \mathbb{R}^d\times\mathcal{M} \rightarrow \mathbb{R}^d$ satisfying
\begin{align}\label{eq:Th2.2-1}
\int_{\mathcal{A}}\left(\frac{p-1}{p \gamma}\right)^{\frac{1}{p-1}}\bigg(H(t, x, i, a, J(t,x,i;\pi))+\sum_{j=1}^L &q_{ij}J(t, x, j;\pi)+\frac{\gamma}{p-1}+\psi_p(t, x, i;\pi)\bigg)_{+}^{\frac{1}{p-1}} \mathrm{d} a=1
\end{align}
where $H(t, x, i, a, J(t, x, i;\pi))$ is defined by \eqref{eq:Hamiltonian}. We consider the following mapping $\mathcal{I}_p$ on $\Pi$ given by, for $\pi \in \Pi$,
\begin{align}\label{eq:Th2.2-2}
\mathcal{I}_p(\pi):=\left(\frac{p-1}{p \gamma}\right)^{\frac{1}{p-1}}\bigg(H(t, x, i, a, J(t, x, i;\pi))+\sum_{j=1}^L &q_{ij}J(t, x, j;\pi)+\frac{\gamma}{p-1}+\psi_p(t, x, i;\pi)\bigg)_{+}^{\frac{1}{p-1}}, p \geq 1
\end{align}
and
$
\mathcal{I}_1(\pi):=\lim _{p \downarrow 1} \mathcal{I}_p(\pi)=\dfrac{\exp \left\{\frac{1}{\gamma} H(t, x, i, a, J(t, x, i;\pi))\right\}}{\int_{\mathcal{A}} \exp \left\{\frac{1}{\gamma} H(t, x, i, a, J(t, x, i;\pi))\right\} \mathrm{d} a}.
$
Denote by $\pi^{\prime}=\mathcal{I}_p(\pi)$ for $\pi \in \Pi$. If $\pi^{\prime} \in \Pi$, then $J\left(t, x, i ; \pi^{\prime}\right) \geq J(t, x, i ; \pi)$ for all $(t, x, i) \in[0, T] \times \mathbb{R}^d \times \mathcal{M}$. Moreover, if the mapping $\mathcal{I}_p: \Pi \rightarrow \Pi$ has a fixed point $\widehat{\pi} \in \Pi$, then $\widehat{\pi}$ is the optimal policy that, for all $(t, x, i) \in[0, T] \times \mathbb{R}^d \times \mathcal{M}$,
$$
V(t, x, i)=\sup_{\pi \in \Pi}J(t, x, i ; \pi)=J\left(t, x, i; \widehat{\pi}\right).
$$
\end{theorem}

\begin{proof}
For two given admissible policies $\pi$ and $\pi^{\prime}$, and any $0 \leq t \leq T$, we consider the discounted process $e^{-\beta s} J(s, \tilde{X}_s^{\pi^{\prime}}, \alpha_s; \pi)$, which captures the value function under the target policy $\pi$ along the behavior trajectory of $\pi^{\prime}$. Now, let the function $H$ defined by \eqref{eq:Hamiltonian} and for $(t, x, i, a) \in[0, T] \times \mathbb{R}^d \times \mathcal{M} \times \mathcal{A}$, set
\begin{align*}
q(t, x, i, a ; \pi):=&H(t, x,i, a, J(t, x, i ; \pi))+\sum_{j=1}^L q_{ij}J(t, x, j; \pi).
\end{align*}
Let $u_n=\inf\{s \geq t:|\tilde{X}_s^{\pi^{\prime}}| \geq n\}$ be a sequence of stopping times. Then, by applying Itô's lemma, we obtain
\begin{align*}
 &e^{-\beta T\wedge u_n} J(T\wedge u_n, \tilde{X}_{T\wedge u_n}^{\pi^{\prime}},\alpha_T; \pi)-e^{-\beta t} J(t, \tilde{X}_t^{\pi^{\prime}},\alpha_t; \pi)\notag\\
 &+\int_t^{T\wedge u_n} e^{-\beta s} \int_{\mathcal{A}}[r(s, \tilde{X}_s^{\pi^{\prime}},\alpha_s, a)+\gamma l_p(\pi^{\prime}(a \mid s, \tilde{X}_s^{\pi^{\prime}}, \alpha_s))]\pi^{\prime}(a \mid s, \tilde{X}_s^{\pi^{\prime}}, \alpha_s) \mathrm{d} a \mathrm{d} s \notag\\
=&\int_t^{T\wedge u_n}e^{-\beta s}\int_{\mathcal{A}}[q(s, \tilde{X}_s^{\pi^{\prime}},\alpha_s,a; \pi)+\gamma l_p(\pi^{\prime}(a \mid s, \tilde{X}_s^{\pi^{\prime}}, \alpha_s))] \pi^{\prime}(a \mid s, \tilde{X}_s^{\pi^{\prime}}, \alpha_s) \mathrm{d} a \mathrm{d} s\notag \\
&+\int_t^{T\wedge u_n} e^{-\beta s}J_{x}(s, \tilde{X}_s^{\pi^{\prime}},\alpha_s; \pi)\circ\tilde{\sigma}(s, \tilde{X}_s^{\pi^{\prime}},\alpha_s, \pi^\prime) \mathrm{d} W_s+\int_t^{T\wedge u_n}e^{-\beta s}\sum_{j=1}^L [J(s, \tilde{X}_s^{\pi^{\prime}},j; \pi)-J(s, \tilde{X}_s^{\pi^{\prime}},\alpha_{s-}; \pi)]\mathrm{d}\widetilde{\Phi}_j(s)\\
\geq &\int_t^{T\wedge u_n}e^{-\beta s}\int_{\mathcal{A}}[q(s, \tilde{X}_s^{\pi^{\prime}},\alpha_s, a; \pi)+\gamma l_p(\pi(a \mid s, \tilde{X}_s^{\pi^{\prime}}, \alpha_s))] \pi(a \mid s, \tilde{X}_s^{\pi^{\prime}}, \alpha_s) \mathrm{d} a \mathrm{d} s\notag \\
&+\int_t^{T\wedge u_n} e^{-\beta s}J_{x}(s, \tilde{X}_s^{\pi^{\prime}},\alpha_s; \pi)\circ\tilde{\sigma}(s, \tilde{X}_s^{\pi^{\prime}},\alpha_s, \pi^\prime) \mathrm{d} W_s+\int_t^{T\wedge u_n}e^{-\beta s}\sum_{j=1}^L [J(s, \tilde{X}_s^{\pi^{\prime}},j; \pi)-J(s, \tilde{X}_s^{\pi^{\prime}},\alpha_{s-}; \pi)]\mathrm{d}\widetilde{\Phi}_j(s)\\
=&\int_t^{T\wedge u_n} e^{-\beta s}J_{x}(s, \tilde{X}_s^{\pi^{\prime}},\alpha_s; \pi)\circ\tilde{\sigma}(s, \tilde{X}_s^{\pi^{\prime}},\alpha_s, \pi^\prime) \mathrm{d} W_s+\int_t^{T\wedge u_n}e^{-\beta s}\sum_{j=1}^L [J(s, \tilde{X}_s^{\pi^{\prime}},j; \pi)-J(s, \tilde{X}_s^{\pi^{\prime}},\alpha_{s-}; \pi)]\mathrm{d}\widetilde{\Phi}_j(s),
\end{align*}
where the inequality and the last equality follow from Lemma \ref{lem:unique-maximizer} and equation \eqref{eq:obj-fun(general)}.
Thus, taking conditional expectation with respect to $\tilde{X}_t^{\pi^{\prime}}=x,\alpha_t=i$ on both sides yields
\begin{align}\label{eq:policy-inprovement-Iteration-proof}
J(t, x, i ; \pi)
\leq & \Ex\bigg[\int_t^{T \wedge u_n} e^{-\beta(s-t)} \int_{\mathcal{A}}[r(s, \tilde{X}_s^{\pi^{\prime}},\alpha_s, a)+\gamma l_p(\pi^{\prime}(a \mid s, \tilde{X}_s^{\pi^{\prime}}, \alpha_s))]\pi^{\prime}(a \mid s, \tilde{X}_s^{\pi^{\prime}}, \alpha_s) \mathrm{d} a \mathrm{d} s \notag\\
&+e^{-\beta\left(T \wedge u_n-t\right)} h(\tilde{X}_{T \wedge u_n}^{\pi^{\prime}}) \mid \tilde{X}_t^{\pi^{\prime}}=x, \alpha_t=i\bigg].
\end{align}

It follows from Assumption \ref{assumption}-(iv), Definition \ref{def:admissible}-(iii) and the moment estimate in Lemma \ref{lem:growth-condition} that there exist constants $\mu^{\prime}, C>0$ such that
\begin{align*}
&\Ex\bigg[\int_t^{T \wedge u_n} e^{-\beta(s-t)} \int_{\mathcal{A}}[r(s, \tilde{X}_s^\pi, \alpha_s, a)+\gamma l_p(\pi^{\prime}(a \mid s, \tilde{X}_s^{\pi^{\prime}}, \alpha_s))]\pi^{\prime}(a \mid s, \tilde{X}_s^{\pi^{\prime}}, \alpha_s) \mathrm{d} a \mathrm{d} s \mid \tilde{X}_t^{\pi'}=x,\alpha_t=i\bigg] \\
\leq &\Ex\bigg[\int_t^T \int_{\mathcal{A}}|r(s, \tilde{X}_s^\pi, \alpha_s, a)+\gamma l_p(\pi^{\prime}(a \mid s, \tilde{X}_s^{\pi^{\prime}}, \alpha_s))|\pi^{\prime}(a \mid s, \tilde{X}_s^{\pi^{\prime}}, \alpha_s)\mathrm{d} a \mathrm{d} s \mid \tilde{X}_t^{\pi'}=x, \alpha_t=i\bigg] \\
\leq &C\int_t^T \Ex\big[(1+|\tilde{X}_s^{\pi'}|^{\mu^{\prime}})\mid \tilde{X}_t^{\pi'}=x, \alpha_t=i\big]\mathrm{d} s
\leq C\int_t^T C(1+|x|^{\mu^{\prime}}) ds
\leq C(1+|x|^{\mu^{\prime}})(T-t)<\infty.
\end{align*}
Additionally, using the moment estimate from Lemma \ref{lem:growth-condition}, we can obtain that there exist constants $\mu^{\prime}, C>0$ such that
\begin{align*}
&\Ex\left[e^{-\beta(T \wedge u_n-t)} h(\tilde{X}_{T \wedge u_n}^{\pi^{\prime}}) \mid \tilde{X}_t^{\pi^{\prime}}=x, \alpha_t=i\right]\\
=&\Ex\left[e^{-\beta (T-t)} h(\tilde{X}_{T}^{\boldsymbol{\pi}^{\prime}})\cdot 1_{\{\max\limits_{t \leq s \leq T}|\tilde{X}_T^{\pi^{\prime}}|<n\}} \mid \tilde{X}_t^{\pi^{\prime}}=x, \alpha_t=i\right]+\Ex\left[e^{-\beta (u_n-t)} h(\tilde{X}_{u_n}^{\pi^{\prime}})\cdot 1_{\{\max\limits_{t \leq s \leq T}|\tilde{X}_T^{\pi^{\prime}}| \geq n\}} \mid \tilde{X}_t^{\pi^{\prime}}=x, \alpha_t=i\right]\\
\leq& \Ex\left[h(\tilde{X}_{T}^{\pi^{\prime}}) \cdot 1_{\{\max\limits_{t \leq s \leq T}|\tilde{X}_T^{\pi^{\prime}}|<n\}} \mid \tilde{X}_t^\pi=x, \tilde{\alpha}_t^{\pi^{\prime}}=i\right]+\Ex\left[ h(\tilde{X}_{u_n}^{\pi^{\prime}}) \cdot 1_{\{\max\limits_{t \leq s \leq T}|\tilde{X}_T^{\pi^{\prime}}| \geq n\}} \mid \tilde{X}_t^\pi=x, \alpha_t=i\right]\\
\leq&  C \Ex\left[(1+|\tilde{X}_T^{\pi^{\prime}}|^{\mu^{\prime}})\cdot 1_{\{\max\limits_{t \leq s \leq T}|\tilde{X}_T^{\pi^{\prime}}|<n\}} \mid \tilde{X}_t^{\pi^{\prime}}=x, \alpha_t=i\right]+ C (1+n^{\mu^{\prime}} ) \Ex\left[1_{\{\max\limits _{t \leq s \leq T}|\tilde{X}_T^{\pi^{\prime}}| \geq n\}} \mid \tilde{X}_t^{\pi^{\prime}}=x, \alpha_t=i\right] \\
\leq&  C \Ex\left[(1+C (1+|x|^{\mu^{\prime}} ))\cdot 1_{\{\max\limits_{t \leq s \leq T}|\tilde{X}_T^{\pi^{\prime}}|<n\}} \mid \tilde{X}_t^{\pi^{\prime}}=x, \alpha_t=i\right]+ C (1+n^{\mu^{\prime}} ) \frac{\Ex\left[\max\limits_{t \leq s \leq T}|\tilde{X}_T^{\pi^{\prime}}|^{\mu^\prime+1} \mid \tilde{X}_t^{\pi^{\prime}}=x, \alpha_t=i\right]}{n^{\mu^\prime+1}}\\
\leq&  C(1+C (1+|x|^{\mu^{\prime}} ))\Ex\left[1_{\{\max\limits_{t \leq s \leq T}|\tilde{X}_T^{\pi^{\prime}}|<n\}} \mid \tilde{X}_t^{\pi^{\prime}}=x, \alpha_t=i\right]+\frac{C (1+n^{\mu^{\prime}} )\cdot C (1+|x|^{\mu^\prime+1} )}{n^{\mu^\prime+1}}\\
\leq&  C(1+C (1+|x|^{\mu^{\prime}} ))\mathbb{P} \left(\max\limits _{t \leq s \leq T} |\tilde{X}_s^\pi |<n \mid \tilde{X}_t^\pi=x, \alpha_t=i \right)+\frac{C(1+n^{\mu^{\prime}})(1+|x|^{\mu^\prime+1})}{n^{\mu^\prime+1}}\\
\leq&  C(1+C(1+|x|^{\mu^{\prime}}))+ \frac{C (1+n^{\mu^{\prime}} ) (1+|x|^{\mu^\prime+1} )}{n^{\mu^\prime+1}}<\infty.
\end{align*}
Thus, applying the dominated convergence theorem and taking the limit as \( n \rightarrow \infty \) on both sides of equation \eqref{eq:policy-inprovement-Iteration-proof}, we obtain: 
\begin{align*}
J(t, x, i ; \pi)
\leq & \Ex\bigg[\int_t^{T} e^{-\beta(s-t)} \int_{\mathcal{A}}\big[r(s, \tilde{X}_s^{\pi^{\prime}},\alpha_s, a)+\gamma l_p(\pi^{\prime}(a \mid s, \tilde{X}_s^{\pi^{\prime}}, \alpha_s))\big]\pi^{\prime}(a \mid s, \tilde{X}_s^{\pi^{\prime}}, \alpha_s) \mathrm{d} a \mathrm{d} s \notag\\
&+e^{-\beta\left(T-t\right)} h(\tilde{X}_{T}^{\pi^{\prime}}) \mid \tilde{X}_t^{\pi^{\prime}}=x, \alpha_t=i\bigg]
= J(t,x,i;\pi').
\end{align*}
This proves that $\pi^{\prime}=\mathcal{I}_p(\pi)$ improves upon $\pi$.

Moreover, if $\mathcal{I}_p(\hat{\pi}) =\hat{\pi}$, then Lemma \ref{lem:unique-maximizer} implies that
\begin{align*}
 &\sup _{\boldsymbol{\pi} \in \mathcal{P}(\mathcal{A})} \int_{\mathcal{A}}\bigg[H(t, x, i, a, J(t,x,i;\hat{\pi}))+\sum_{j=1}^L q_{ij}J(t, x, j;\hat{\pi})+\gamma l_p(\pi(a))\bigg] \pi(a) \mathrm{d} a\\
 &=\int_{\mathcal{A}}\bigg[H(t, x, i, a, J(t, x,i;\hat{\pi}))+\sum_{j=1}^L q_{ij}J(t, x, j;\hat{\pi})+\gamma l_p( \hat{\pi}(a \mid t, x, i))\bigg] \hat{\pi}(a \mid t, x, i) \mathrm{d} a =0,
\end{align*}
where the last equality follows from \eqref{eq:obj-fun(general)} with $\pi(a \mid t, x, i)$ replaced by $\hat{\pi}(a \mid t, x, i)$.
This means that $J(t,x,i;\hat{\pi})$ satisfies the exploratory HJB equation \eqref{eq:obj-fun(general)} and so $J(t,x,i;\hat{\pi})$ is the optimal value function and hence $\hat{\pi}$ is the optimal policy.
\end{proof}
\begin{remark}
Theorem \ref{thm:policy-improve} is a theoretical result that cannot be directly applied to learning procedures, because the Hamiltonian $H(t, x, i, a, J(t, x, i;\pi))$ relies on the knowledge of the model parameters, which we do not have in the RL context. Consequently, we propose a model-free RL algorithm by extending the q-learning theory of \citet{jia2023q} to fit our formulation under Tsallis entropy.
\end{remark}

\section{q-function for Markov Regime-switching System and Martingale Characterization under Tsallis Entropy}


Following the framework of \citet{jia2023q}, we define a Q-function parameterized by a time step $\Delta t > 0$. For any policy $\pi \in \Pi$ and $(t, x, i, a) \in[0, T) \times \mathbb{R}^d\times\mathcal{M} \times \mathcal{A}$, we construct a perturbed policy $\tilde{\pi}$ that enforces the constant action $a$ on $[t, t+\Delta t)$ and then follows $\pi$ on $[t+\Delta t, T]$. The corresponding state process $X^{\tilde{\pi}}$, initialized at $X_t^{\tilde{\pi}}=x$ is piece wise continuous. Specifically, on $[t, t+\Delta t)$, the trajectory $X^a$ satisfies the stochastic differential equation
\begin{align}\label{eq:X-a}
\mathrm{d} X_s^a=b(s, X_s^a, \alpha_s, a) \mathrm{d} s+\sigma(s, X_s^a,\alpha_s, a) \mathrm{d} W_s,\; s \in[t, t+\Delta t),\; X_t^a=x,\alpha_s=i.
\end{align}
On $[t+\Delta t, T]$, the process evolves according to \eqref{eq:state-Exp} with the initial condition $(t+\Delta t, X_{t+\Delta t}^a, \alpha_{t+\Delta t})$.
The $\Delta t$-parameterized Q-function for policy $\pi$ is defined as the expected total reward from the perturbed policy $\tilde{\pi}$:
\begin{align}\label{eq:Q-delta-func}
 Q_{\Delta t}(t, x,i, a ; \pi)
=&\Ex\bigg[\int_t^{t+\Delta t} e^{-\beta(s-t)} r(s, X_s^a,\alpha_s, a) \mathrm{d} s\notag\\
&+\int_{t+\Delta t}^T e^{-\beta(s-t)}r(s, X_s^{\pi},\alpha_s, a_s^{\pi}) \mathrm{d} s+e^{-\beta(T-t)}h(X_T^{\pi}) \mid X_t^{\tilde{\pi}}=x,\alpha_t=i\bigg].
\end{align}
The following proposition provides an asymptotic expansion of this Q-function with respect to $\Delta t$.
\begin{proposition}\label{pro:Q-func}
We have
\begin{align}\label{eq:pro-Q-func}
 Q_{\Delta t}(t, x,i, a; \pi)
=  H(t, x,i, a, J(t, x,i ; \pi))+\sum_{j=1}^L q_{ij}J(t, x, j; \pi)] \Delta t+o(\Delta t) .
\end{align}
\end{proposition}
\begin{proof}
Applying Itô's lemma, we obtain
\begin{align*}
    &\Ex\bigg[e^{-\beta \Delta t}J(t+\Delta t, X_{t+\Delta t }^a,\alpha_{t+\Delta t}; \pi)-J(t, x,i ; \pi)\mid X_t^{\tilde{\pi}}=x,\alpha_t=i\bigg]\\
    =&\Ex\bigg[\int_t^{t+\Delta t}e^{-\beta(s-t)}[H(s, X_s^a,\alpha_s, a, J(s, X_s^a,\alpha_s; \pi))+\sum_{j=1}^L q_{ij}J(s, X_s^{a}, j; \pi)] \mathrm{d} s\mid X_t^{\tilde{\pi}}=x,\alpha_t=i\bigg].
\end{align*}

Next, breaking the time period $[t, T]$ into $[t, t+\Delta t)$ and $[t+\Delta t, T]$, and conditioning on the state at $t+\Delta t$ for the period $[t+\Delta t, T]$ lead to:
\begin{align*}
 &Q_{\Delta t}(t, x,i, a ; \pi) \\
= &\Ex\Bigg[\int_t^{t+\Delta t} e^{-\beta(s-t)} r(s, X_s^a,\alpha_s, a) \mathrm{d} s+\Ex\bigg[\int_{t+\Delta t}^T e^{-\beta(s-t)}(r(s, X_s^{\pi},\alpha_s, a_s^{\pi})\\
&\quad+\gamma l_p(\pi(a\mid s, X_s^a,\alpha_s)) )\mathrm{d} s+e^{-\beta(T-t)}h(X_T^{\pi}) \mid X_{t+\Delta t}^a, \alpha_{t+\Delta t}\bigg] \mid X_t^{\tilde{\pi}}=x,\alpha_t=i\Bigg] \\
= & \Ex\bigg[\int_t^{t+\Delta t} e^{-\beta(s-t)} r(s, X_s^a,\alpha_s, a) \mathrm{d} s+e^{-\beta \Delta t} J(t+\Delta t, X_{t+\Delta t}^a,\alpha_{t+\Delta t};\pi)\mid X_t^{\tilde{\pi}}=x,\alpha_t=i\bigg]\\
= & \Ex\bigg[\int_t^{t+\Delta t} e^{-\beta(s-t)} r(s, X_s^a,\alpha_s, a) \mathrm{d} s+e^{-\beta \Delta t} J(t+\Delta t, X_{t+\Delta t}^a,\alpha_{t+\Delta t}; \pi) -J(t, x, i;\pi)\mid X_t^{\tilde{\pi}}=x,\alpha_t=i\bigg]+J(t, x, i;\pi)\\
=&\Ex\Bigg[\int_{t}^{t+\Delta t}e^{-\beta(s-t)}\bigg[H(s, X_s^a,\alpha_s, a, J(s, X_s^a,\alpha_s; \pi)) +\sum_{j=1}^L q_{\alpha_sj}J(s, X_s^{a}, j; \pi)\bigg]\mathrm{d} s \mid X_t^{\tilde{\pi}}=x,\alpha_t=i\Bigg]+J(t, x, i; \pi) \\
=&J(t, x, i ; \pi)+\bigg[H(t, x, i, a, J(t, x, i; \pi))+\sum_{j=1}^L q_{ij}J(t, x, j; \pi)\bigg]\Delta t+o(\Delta t),
\end{align*}
where the second-to-last equality follows from It{\^o}'s lemma and the error order is due to the approximation of the integral involved.
\end{proof}

Next, we define the q-function as the counterpart of the Q-function in the continuous-time framework and provide its martingale characterization.

\begin{definition}[q-function]\label{def:q-corollary}
 The q-function of the problem \eqref{eq:state-Exp}-\eqref{eq:value-fun(MC)} associated with a given policy $\pi \in \Pi$ is defined as
\begin{align}\label{eq:q-def}
q(t, x, i, a ; \pi)=&H(t, x, i, a, J(t, x, i ; \pi))+\sum_{j=1}^L q_{ij}J(t, x, j; \pi),\quad
(t, x, i, a) \in[0, T] \times \mathbb{R}^d \times \mathcal{M} \times \mathcal{A}.
\end{align}
\end{definition}
As an immediate consequence of Proposition \ref{pro:Q-func}, the q-function can be characterized as the first-order derivative of the Q-function with respect to $\Delta t$, that is:
\begin{align}\label{eq:q-corollary}
q(t, x,i, a ; \pi)=\lim _{\Delta t \rightarrow 0} \frac{Q_{\Delta t}(t, x,i, a ; \pi)-J(t, x ,i; \pi)}{\Delta t} .
\end{align}

\begin{remark}
The improved policy $\pi^{\prime}$ in Theorem \ref{thm:policy-improve} can be expressed in terms of the q-function as:
$$
\pi^{\prime}(a \mid t, x, i)=\left(\frac{p-1}{p \gamma}\right)^{\frac{1}{p-1}}\big(q(t, x, i, a ; \pi)+\frac{\gamma}{p-1}+\psi_p(t, x, i;\pi)\big)_{+}^{\frac{1}{p-1}}, \quad \forall p \geq 1,
$$
where the Lagrange multiplier $\psi_p(t, x, i;\pi):[0, T] \times \mathbb{R}^d \times\mathcal{M} \mapsto \mathbb{R}$ must satisfy the normalization condition:
\begin{align}\label{eq:Remark3.1-1}
\int_{\mathcal{A}}\left(\frac{p-1}{p \gamma}\right)^{\frac{1}{p-1}}\big(q(t, x,i, a ;\pi)+\frac{\gamma}{p-1}+\psi_p(t, x, i;\pi)\big)_{+}^{\frac{1}{p-1}} \mathrm{d} a=1 .
\end{align}
To address the natural question of whether such a function $\psi(t, x, i;\pi)$ exists, we consider the map $u \mapsto F(u)$ for a given policy $\pi$ and a fixed triplet $(t, x, i) \in [0, T] \times \mathbb{R}^d\times\mathcal{M}$:
$$
u \mapsto F(u):=\int_{\mathcal{A}}\left(\frac{p-1}{p \gamma}\right)^{\frac{1}{p-1}}(q(t, x, i, a ; \pi)+\frac{\gamma}{p-1}+u)_{+}^{\frac{1}{p-1}} \mathrm{d} a .
$$
Provided the q-function satisfies necessary integral conditions for $F(u)$ to be well-defined, the mapping $u \mapsto F(u)$ is continuous and strictly increasing, with $F(u) \rightarrow 0$ as $u \rightarrow-\infty$ and $F(u) \rightarrow+\infty$ as $u \rightarrow+\infty$. This guarantees the existence and uniqueness of the function $\psi_p(t, x, i;\pi)$ satisfying \eqref{eq:Remark3.1-1}.
\end{remark}

The following result provides a martingale characterization of the q-function under a given policy $\pi$ when the value function is specified.
\begin{proposition}\label{pro:martingale}
Let $\pi\in\Pi$ be a policy whose value function $J(\cdot, \cdot, i; \pi) \in C^{1,2}([0, T) \times \mathbb{R}^d) \cap C([0, T] \times \mathbb{R}^d)$ for all $i\in \mathcal{M}$ satisfies \eqref{eq:obj-fun(general)}. Let a continuous function $\hat{q}:[0, T] \times \mathbb{R}^d\times \mathcal{M} \times \mathcal{A}\mapsto \mathbb{R}$ be given. Then, we have
\begin{itemize}
\item[(i)]$\hat{q}(t, x, i, a)=q(t, x, i, a ; \pi)$ for all $(t, x, i, a) \in [0, T] \times \mathbb{R}^d \times\mathcal{M} \times \mathcal{A}$ if and only if for all $(t, x, i) \in[0, T] \times \mathbb{R}^d \times\mathcal{M}$, the following process
\begin{align}\label{eq:Pro3.2-1}
e^{-\beta s}J(s, X_s^{\pi},\alpha_s ; \pi)+\int_t^s e^{-\beta u}[r(u, X_u^{\pi},\alpha_u, a_u^{\pi})-\hat{q}(u, X_u^{\pi},\alpha_u, a_u^{\pi})] \mathrm{d} u,
\end{align}
is an $(\mathbb{F}, \mathbb{P})$-martingale, where $\{X_s^{\pi},\; t\leq s\leq T\}$ is the solution to \eqref{eq:state-Exp} under $\pi$ with $X_t^{\pi}=x,\; \alpha_t=i$.
\item[(ii)] If $\hat{q}(t, x, i, a)=q(t, x,  i, a; \pi)$ for all $(t, x,  i, a) \in[0, T] \times \mathbb{R}^d \times\mathcal{M} \times \mathcal{A}$, given any $\pi^{\prime} \in \Pi$, for all $(t, x, i) \in[0, T] \times \mathbb{R}^d \times\mathcal{M}$, the following process
\begin{align}\label{eq:Pro3.2-2}
e^{-\beta s}J(s, X_s^{\pi^{\prime}}, \alpha_s; \pi)+\int_t^s e^{-\beta u}[r(u, X_u^{\pi^{\prime}}, \alpha_u, a_u^{\pi^{\prime}})-\hat{q}(u, X_u^{\pi^{\prime}}, \alpha_u, a_u^{\pi^{\prime}})]\mathrm{d} u,
\end{align}
is an $(\mathbb{F}, \mathbb{P})$-martingale. Here, $\{X_s^{\pi_p^{\prime}},t\leq s\leq T\}$ is the solution to \eqref{eq:state-Exp} under $\pi^{\prime}$ with $X_t^{\pi^{\prime}}=x,\; \alpha_t=i$.
\item[(iii)] If there exists $\pi^{\prime} \in \Pi$ such that for all $(t, x, i) \in[0, T] \times \mathbb{R}^d \times\mathcal{M}$, the process \eqref{eq:Pro3.2-2} is an $(\mathbb{F}, \mathbb{P})$-martingale with initial condition $X_t^{\pi^{\prime}}=x, \alpha_t=i$, then $\hat{q}(t, x, i, a)= q(t, x, i, a; \pi)$ for all $(t, x, i, a) \in[0, T] \times\mathbb{R}^d \times\mathcal{M} \times \mathcal{A}$.
\end{itemize}
Moreover, in any of the three cases above, the q-function satisfies
\begin{align}\label{eq:q-fun-satisfies}
\int_{\mathcal{A}}[q(t, x, i, a; \pi)+\gamma l_p(\pi(a \mid t, x, i))] \pi(a \mid t, x, i) \mathrm{d} a=0.
\end{align}
\end{proposition}

\begin{proof}
(i) To simplify our notation, we first set
\begin{align*}
    M(s,Y, Z;\pi):=&\int_t^s e^{-\beta u}\sum_{j=1}^{L}[J(u, Y_u,j; \pi)-J(u, Y_u,\alpha_{u-}; \pi)]\mathrm{d}\widetilde{\Phi}_j(u)\\
    &\quad +\int_t^s e^{-\beta u} J_{x} (u, Y_u, \alpha_u ; \pi) \circ \sigma(u, Y_u, \alpha_u, Z_u) \mathrm{d} W_u.
\end{align*}
Applying Itô's lemma to $e^{-\beta s} J(s, X_s^{\pi}, \alpha_s ; \pi)$, we obtain for $0 \leq t<s \leq T$ :
\begin{align}\label{eq:35}
& e^{-\beta s} J(s, X_s^{\pi}, \alpha_s; \pi)+\int_t^s e^{-\beta u}[r(u, X_u^{\pi}, \alpha_u, a_u^{\pi})-\hat{q}(u, X_u^{\pi}, \alpha_u, a_u^{\pi})] \mathrm{d} u \notag\\
=&e^{-\beta t} J(t, X_t^{\pi}, \alpha_t; \pi)+\int_t^s e^{-\beta u}\bigg[H(u, X_u^{\pi}, \alpha_u, a_u^{\pi}, J(u, X_u^{\pi}, \alpha_u ; \pi))+\sum_{j=1}^{L} q_{ij}J(u, X_u^{\pi}, j ; \pi) \notag\\
&-\hat{q}(u, X_u^{\pi}, \alpha_u, a_u^{\pi})\bigg] \mathrm{d} u + M(s,X^\pi,a^\pi;\pi)\notag\\
=&e^{-\beta t} J(t, X_t^{\pi}, \alpha_t; \pi)+\int_t^s e^{-\beta u}\bigg[q(u, X_u^{\pi}, \alpha_u, a_u^{\pi} ; \pi)-\hat{q}(u, X_u^{\pi}, \alpha_u, a_u^{\pi})\bigg]\mathrm{d} u+ M(s,X^\pi,a^\pi;\pi).
\end{align}
Thus, if $\hat{q} \equiv q(\cdot, \cdot, \cdot, \cdot ; \pi)$, then we have 

\begin{align*}
e^{-\beta s} J(s, X_s^{\pi}, \alpha_s ; \pi)+\int_t^s e^{-\beta u}[r(u, X_u^{\pi}, \alpha_u, a_u^{\pi})-\hat{q}(u, X_u^{\pi}, \alpha_u, a_u^{\pi})] \mathrm{d} u
=e^{-\beta t} J(t, X_t^{\pi}, \alpha_t; \pi)+ M(s,X^\pi,a^\pi;\pi).
\end{align*}
Therefore, \eqref{eq:Pro3.2-1} is an $(\mathbb{F}, \mathbb{P})$-martingale on $[t, T]$.

Conversely, if the left hand side of \eqref{eq:35} is a martingale, then by observing that $M(s,X^\pi,a^\pi;\pi)$ is an $(\mathbb{F}, \mathbb{P})$-martingale, it follows directly from \eqref{eq:35} that 
$\int_t^s e^{-\beta u}[q(u, X_u^\pi, \alpha_u, a_u^\pi; \pi)-\hat{q}(u, X_u^\pi, \alpha_u, a_u^\pi)] \mathrm{d} u$ is a continuous martingale with finite variation and hence zero quadratic variation. Therefore,  $\int_t^s e^{-\beta u}[q(u, X_u^{\boldsymbol{\pi}}, \alpha_u, a_u^{\boldsymbol{\pi}} ; \pi)-\hat{q}(u, X_u^{\boldsymbol{\pi}}, \alpha_u, a_u^{\boldsymbol{\pi}})] \mathrm{d} u=0$ $\mathbb{P}$-almost surely for all $s \in[t, T]$; see, e.g., (\citet{Karatzas1991}, Chapter 1, Exercise 5.21)


Denote $f(t, x, i, a)=q(t, x, i, a ; \pi)-\hat{q}(t, x, i, a)$. Then $f$ is a continuous function that maps $[0, T] \times \mathbb{R}^d \times \mathcal{M}\times \mathcal{A}$ to $\mathbb{R}$. Suppose the desired conclusion is not true, then there exists a quadruple $( t^*, x^*, i^*, a^*)$ and $\epsilon>0$ such that $f(t^*, x^*, i^*, a^*)>\epsilon$. Because $f$ is continuous, there exists $\delta>0$ such that $f(u, x^{\prime}, i^*, a^{\prime})>\epsilon/2$ for all $(u, x^{\prime}, i^*, a^{\prime})$ with $|u-t^*| \vee|x^{\prime}-x^*|\vee|a^{\prime}-a^*|<\delta$. Here “$\vee$” means taking the larger one, i.e., $u \vee v=\max \{u, v\}$. Now consider the state process, still denoted by $X^\pi$, starting from $(t^*, x^*, i^*, a^*)$, namely, $\{X_s^\pi, t^* \leq s \leq T\}$ follows \eqref{eq:state-Exp} with $X_{t^*}^\pi=x^*, \alpha_{t^*}=i^*$ and $a_{t *}^\pi=a^*$. Define
$$
\tau=\inf \left\{u \geq t^*:|u-t^*|>\delta \text{ or }|X_u^\pi-x^*|>\delta  \text{ or }  \alpha_{u}\neq i^*\right\}.
$$
The continuity of $X^\pi$ implies that $u>t^*$, $\mathbb{P}$-almost surely. Here “$\wedge$” means taking the smaller one, i.e., $u \wedge v=\min \{u, v\}$. We have already proved that there exists $\Omega_0 \in \mathcal{F}$ with $\mathbb{P}(\Omega_0)=0$ such that for all $\omega \in \Omega \backslash \Omega_0$, $\int_{t^*}^s e^{-\beta u} f(u, X_u^\pi(\omega),\alpha_u(\omega), a_u^\pi(\omega)) \mathrm{d} u=0$ for all $s \in[t^*, T]$. It follows from Lebesgue's differentiation theorem that for any $\omega \in \Omega \backslash \Omega_0$,
$$
f(s, X_s^\pi(\omega), \alpha_s(\omega), a_s^\pi(\omega))=0, \text { a.e. } s \in[t^*, \tau(\omega)].
$$

Consider the set $Z(\omega)=\{s \in[t^*, \tau(\omega)]: a_s^\pi(\omega) \in \mathcal{B}_\delta(a^*)\} \subset[t^*, \tau(\omega)]$, where $\mathcal{B}_\delta(a^*)=\{a^{\prime} \in \mathcal{A}$ : $|a^{\prime}-a^*|>\delta\}$ is the neighborhood of $a^*$. Because $f(s, X_s^\pi(\omega), \alpha_s(\omega), a_s^\pi(\omega))>\frac{\epsilon}{2}$ when $s \in Z(\omega)$, we conclude that $Z(\omega)$ has Lebesgue measure zero for any $\omega \in \Omega \backslash \Omega_0$. That is,
$$
\int_{[t^*, T]} \mathbf{1}_{\{s \leq \tau(\omega)\}} \mathbf{1}_{\{a_x^\pi(\omega) \in \mathcal{B}_\delta(a^*)\}} \mathrm{d} s=0.
$$
Integrating $\omega$ with respect to $\mathbb{P}$ and applying Fubini's theorem, we obtain
\begin{align*}
0 & =\int_{\Omega} \int_{t^*}^T \mathbf{1}_{\{s \leq \tau(\omega)\}} \mathbf{1}_{\{a_x^\pi(\omega) \in \mathcal{B}_\delta(a^*)\}} \mathrm{d} s \mathbb{P}(\mathrm{d} \omega)=\int_{t^*}^T\int_{\Omega}  \mathbf{1}_{\{s \leq \tau(\omega)\}} \mathbf{1}_{\{a_x^\pi(\omega) \in \mathcal{B}_\delta(a^*)\}}\mathbb{P}(\mathrm{d} \omega) \mathrm{d} s \\
&=\int_{t^*}^T \Ex[\mathbf{1}_{\{s \leq \tau\}} \mathbf{1}_{\{a_x^\pi \in \mathcal{B}_\delta(a^*)\}}] \mathrm{d} s =\int_{t^*}^T \Ex\left[\Ex\left[\mathbf{1}_{\{s \leq \tau\}} \mathbf{1}_{\{a_x^\pi \in \mathcal{B}_\delta(a^*)\}}\mid \mathcal{F}_s\right]  \right] \mathrm{d} s\\
&=\int_{t^*}^T \Ex\left[\mathbf{1}_{\{s \leq \tau\}}\cdot \Ex\left[\mathbf{1}_{\{a_x^\pi \in \mathcal{B}_\delta(a^*)\}}\mid \mathcal{F}_s\right]  \right] \mathrm{d} s
 =\int_{t^*}^T \Ex\left[\mathbf{1}_{\{s \leq \tau\}} \mathbb{P}(a_s^\pi \in \mathcal{B}_\delta(a^*) \mid \mathcal{F}_s)\right] \mathrm{d} s\\
&=\int_{t^*}^T \Ex\left[\mathbf{1}_{\{s \leq \tau\}} \int_{\mathcal{B}_\delta(a^*)} \pi(a \mid s, X_s^\pi,\alpha_s) \mathrm{d} a\right] \mathrm{d} s\\
&\geq \min_{|u-t^*|<\delta ,|x^{\prime}-x^*|<\delta}\left\{\int_{\mathcal{B}_\delta\left(a^*\right)} \pi\left(a \mid u, x^{\prime}, i^*\right) \mathrm{d} a\right\} \int_{t^*}^T \Ex\left[\mathbf{1}_{\{s \leq \tau\}}\right] \mathrm{d} s \\
&= \min_{|u-t^*|<\delta ,|x^{\prime}-x^*|<\delta}\left\{\int_{\mathcal{B}_\delta\left(a^*\right)} \pi\left(a \mid u, x^{\prime}, i^*\right) \mathrm{d} a\right\} \Ex\left[\int_{t^*}^T \mathbf{1}_{\{s \leq \tau\}} \mathrm{d} s\right]\\
& =\min_{|u-t^*|<\delta ,|x^{\prime}-x^*|<\delta}\left\{\int_{\mathcal{B}_\delta(a^*)} \pi(a \mid u, x^{\prime}, i^*) \mathrm{d} a\right\} \Ex\left[(\tau \wedge T)-t^*\right] \geq 0.
\end{align*}

Since $\tau>t^*$, $\mathbb{P}$-almost surely, the above implies
$\min_{|u-t^*|<\delta ,|x^{\prime}-x^*|<\delta}\left\{\int_{\mathcal{B}_\delta\left(a^*\right)} \pi(a \mid u, x^{\prime}, i^*) \mathrm{d} a\right\}=0.$
However, this contradicts Definition \ref{def:admissible} about an admissible policy. Indeed, Definition \ref{def:admissible}-(i) stipulates $\operatorname{supp} \pi(\cdot \mid t, x, i)=\mathcal{A}$ for any $(t, x, i)\in[0,T]\times \mathbb{R}^d \times \mathcal{M}$; hence $\int_{\mathcal{B}_\delta\left(a^*\right)} \pi(a \mid t, x, i) \mathrm{d} a>0$. Then the continuity in Definition \ref{def:admissible}-(iii) yields $
\min_{|u-t^*|<\delta ,|x^{\prime}-x^*|<\delta}\left\{\int_{\mathcal{B}_\delta(a^*)} \pi\left(a \mid u, x^{\prime}, i^*\right) \mathrm{d} a\right\}>0,$
a contradiction. Hence we conclude $q(t, x, a, i ; \pi)=\hat{q}(t, x, a, i)$ for every $(t, x, a, i)$.

(ii) Applying Itô's lemma to $e^{-\beta s} J(s, X_s^{\pi^{\prime}}, \alpha_s; \pi)$, we get
\begin{align*}
& e^{-\beta s} J(s, X_s^{\pi^{\prime}}, \alpha_s ; \pi)+\int_t^s e^{-\beta u}[r(u, X_u^{\pi^{\prime}}, \alpha_u, a_u^{\pi^{\prime}})-\hat{q}(u, X_u^{\pi^{\prime}}, \alpha_u,
 a_u^{\pi^{\prime}})] \mathrm{d} u \\
= &e^{-\beta t} J(t, x, i; \pi)+\int_t^s e^{-\beta u}[q(u, X_u^{\pi^{\prime}}, \alpha_u, a_u^{\pi^{\prime}} ; \pi)-\hat{q}(u, X_u^{\pi^{\prime}}, \alpha_u, a_u^{\pi^{\prime}})] \mathrm{d} u+M(s,X^{\pi'},a^{\pi'};\pi)
\end{align*}
So, when $\hat{q} \equiv q(\cdot, \cdot, \cdot, \cdot; \boldsymbol{\pi})$, the left-hand side of the above equation, i.e., the process defined in \eqref{eq:Pro3.2-2}, is an $(\mathbb{F}, \mathbb{P})$-martingale on $[t, T]$.

(iii) Let $\boldsymbol{\pi}^{\prime} \in \boldsymbol{\Pi}$ be given satisfying the assumption. 
It then follows from (ii) that $\int_t^s e^{-\beta u}[\hat{q}(u, X_u^{\pi^{\prime}}, \alpha_u, a_u^{\pi^{\prime}})-q(u, X_u^{\pi^{\prime}}, \alpha_u, a_u^{\pi^\prime} ; \pi)] \mathrm{d} u$ is an $(\mathbb{F}, \mathbb{P})$-martingale. If the desired conclusion is not true, then the same argument in (i) still applies to conclude that
$\min_{|x^\prime-x^*|<\delta,|u-t^*|<\delta}\left\{\int_{\mathcal{B}_\delta(a^*)} \pi^\prime(a \mid u, x^\prime, i^*) \mathrm{d} a\right\}=0, $
which is a contradiction because $\pi^\prime$ is admissible.
\end{proof}

Similar to Theorem 2 in \citet{jia2023q}, we can strengthen Proposition \ref{pro:martingale} and characterize the q-function and the value function associated with a given policy $\pi$ simultaneously.

\begin{theorem}\label{thm:martingale-characterization}
For each $p \geq 1$, let a policy $\pi \in \Pi$, a function $\hat{J}(\cdot, \cdot, i) \in C^{1,2}([0, T) \times \mathbb{R}^d) \cap C([0, T] \times \mathbb{R}^d)$ for all $i\in \mathcal{M}$ and a continuous function $\hat{q}:[0, T] \times \mathbb{R}^d \times\mathcal{M}\times \mathcal{A} \mapsto \mathbb{R}$ be given such that, for all $(t, x, i, a) \in[0, T] \times \mathbb{R}^d \times \mathcal{M}\times \mathcal{A}$,
\begin{align}\label{eq:The3.3-1}
\hat{J}(T, x, i)=h(x), \quad \int_{\mathcal{A}}\left[\hat{q}(t, x, i, a)+\gamma l_p(\pi(a \mid t, x, i))\right] \pi(a \mid t, x, i) \mathrm{d} a=0 .
\end{align}
Then,
\begin{itemize}
\item[(i)] $\hat{J}(t,x,i)$ and $\hat{q}(t,x,i,a)$ are respectively the value function satisfying Eq.\eqref{eq:obj-fun(general)} and the $q$-function associated with $\pi$ if and only if for all $(t, x, i) \in[0, T] \times \mathbb{R}^d \times\mathcal{M}$, the following process
\begin{align}\label{eq:The3.5-1}
e^{-\beta s}\hat{J}(s, X_s^{\pi},\alpha_s)+\int_t^s e^{-\beta u}[r(u, X_u^{\pi},\alpha_u, a_u^{\pi})-\hat{q}(u, X_u^{\pi},\alpha_u, a_u^{\pi})] \mathrm{d} u,
\end{align}
is an $(\mathbb{F}, \mathbb{P})$-martingale, where $\{X_s^{\pi},t\leq s\leq T\}$ is the solution to \eqref{eq:state-Exp} under $\pi$ with $X_t^{\pi}=x,\; \alpha_t=i$.
\item[(ii)] If $\hat{J}(t,x,i)$ and $\hat{q}(t,x,i,a)$ are respectively the value function and the $q$-function associated with $\pi$, then for any $\pi^{\prime} \in \Pi$ and all $(t, x, i) \in[0, T] \times \mathbb{R}^d \times\mathcal{M}$, the following process
\begin{align}\label{eq:The3.5-2}
e^{-\beta s}\hat{J}(s, X_s^{\pi^{\prime}},\alpha_s)+\int_t^s e^{-\beta u}[r(u, X_u^{\pi^{\prime}},\alpha_u, a_u^{\pi^{\prime}})-\hat{q}(u, X_u^{\pi^{\prime}},\alpha_u, a_u^{\pi^{\prime}})] \mathrm{d} u,
\end{align}
is an $(\mathbb{F}, \mathbb{P})$-martingale, where $\{X_s^{\pi^{\prime}},t\leq s\leq T\}$ is the solution to \eqref{eq:state-Exp} under $\pi^{\prime}$ with $X_t^{\pi^{\prime}}=x,\; \alpha_t=i$.
\item[(iii)] If there exists $\pi^{\prime} \in \Pi$ such that for all $(t, x, i) \in[0, T] \times \mathbb{R}^d \times\mathcal{M}$, the process \eqref{eq:The3.5-2} is an $(\mathbb{F}, \mathbb{P})$-martingale with initial condition $X_t^{\pi^{\prime}}=x, \alpha_t=i$, then $\hat{J}(t,x,i)$ and $\hat{q}(t,x,i,a)$  are respectively the value function and the $q$-function associated with $\pi$.
\end{itemize}
Moreover, in any of the above three cases, if it holds further that
$$
\pi(a \mid t, x, i)=\left(\frac{p-1}{p \gamma}\right)^{\frac{1}{p-1}}(\hat{q}(t, x, i, a)+\frac{\gamma}{p-1}+\psi_p(t, x,i))_{+}^{\frac{1}{p-1}}, \quad p\geq1,
$$
with the normalizing function $\psi_p(t, x,i)$ satisfying $\int_{\mathcal{A}} \left(\frac{p-1}{p \gamma}\right)^{\frac{1}{p-1}}(\hat{q}(t, x, i, a)+\frac{\gamma}{p-1}+\psi_p(t, x,i))_{+}^{\frac{1}{p-1}} \mathrm{d} a=1$ for all $(t, x, i, a) \in[0, T] \times \mathbb{R}^d \times\mathcal{M}\times \mathcal{A}$, then $\pi$ for each $p\geq1$ is an optimal policy and $\hat{J}$ is the corresponding optimal value function.
\end{theorem}

\begin{proof}
(i) The “$\Rightarrow$” part follows directly from the same argument as in the proof of \cref{pro:martingale}. We only need to prove the “$\Leftarrow$” part.
If $e^{-\beta s}\hat{J}(s, X_s^{\pi},\alpha_s ; \pi)+\int_t^s e^{-\beta u}(r(u, X_u^{\pi},\alpha_u, a_u^{\pi})-\hat{q}(u, X_u^{\pi},\alpha_u, a_u^{\pi})) \mathrm{d} u$ is an $(\mathbb{F}, \mathbb{P})$-martingale. Hence, for any initial state $(t, x, i)$, we have
\begin{align*}
\Ex\bigg[e^{-\beta T} \hat{J}(T, X_T^{\pi},\alpha_T)+\int_t^T e^{-\beta s}(r(s, X_s^{\pi}, \alpha_s, a_s^{\pi})-\hat{q}(s, X_s^{\pi},\alpha_s, a_s^{\pi}))\mathrm{d} s \mid \mathcal{F}_t\bigg]=e^{-\beta t} \hat{J}(t, x, i).
\end{align*}
We integrate over the action randomization with respect to the policy $\pi$, and then obtain
\begin{align*}
\Ex\bigg[e^{-\beta T} \hat{J}(T, \tilde{X}_T^{\pi},\alpha_T)+\int_t^T \int_{\mathcal{A}}e^{-\beta s}\big[r(s, \tilde{X}_s^{\pi}, \alpha_s, a)-\hat{q}(s, \tilde{X}_s^{\pi},\alpha_s, a))\pi(a|s,\tilde{X}_s^{\pi},\alpha_s)\big]\mathrm{d}a\mathrm{d} s \mid \mathcal{F}_t^W\bigg]=e^{-\beta t} \hat{J}(t, x, i).
\end{align*}

This, together with the terminal condition $\hat{J}(T, x, i)=h(x)$, and constraint \eqref{eq:The3.3-1} yields that
\begin{align*}
&\hat{J}(t, x, i)
=\Ex\bigg[e^{-\beta(T-t)} h(\tilde{X}_T^{\pi})+\int_t^T e^{-\beta(s-t)} \int_{\mathcal{A}}\big[r(s, \tilde{X}_s^{\pi}, \alpha_s, a)+\gamma l_p(\pi(a \mid s, \tilde{X}_s^{\pi}, \alpha_s))\big]\pi(a \mid s, \tilde{X}_s^{\pi}, \alpha_s) \mathrm{d} a \mathrm{d} s \mid \mathcal{F}_t^W\bigg].
\end{align*}
Hence $\hat{J}(t, x, i)=J(t, x, i; \pi)$ for all $(t, x, i) \in[0, T] \times \mathbb{R}^d \times \mathcal{M}$. Furthermore, based on \cref{pro:martingale}, the martingale condition implies that $\hat{q}(t, x, i, a)=q(t, x, i, a ; \pi)$ for all $(t, x, i, a) \in[0, T] \times \mathbb{R}^d \times\mathcal{M}\times \mathcal{A}$.

(ii) This follows immediately from \cref{pro:martingale}-(ii).

(iii) Let $\pi^{\prime} \in \Pi$ be given satisfying the assumption of this part and set
\begin{align*}
\hat{r}(t, x, i, a):=\hat{J}_t(t, x, i)+b(t, x, i, a) \circ \hat{J}_x(t, x, i)+\frac{1}{2} \sigma \sigma^{\top}(t, x, i, a) \circ \hat{J}_{xx}(t, x, i)
-\beta \hat{J}(t, x, i)+\sum_{j=1}^L q_{i j} \hat{J}(t, x, j).
\end{align*}
Applying Itô's lemma to $e^{-\beta s} \hat{J}(s, X_s^{\pi^{\prime}}, \alpha_s)$, it is easy to verify that $e^{-\beta s} \hat{J}(s, X_s^{\pi^{\prime}}, \alpha_s)-\int_t^s e^{-\beta u}\hat{r}(u, X_u^{\pi^{\prime}}, \alpha_u, a_u^{\pi^{\prime}}) \mathrm{d} u$
is an $(\mathbb{F}, \mathbb{P})$-martingale.
Thus, combining the assumption that the process \eqref{eq:The3.5-2} is an $(\mathbb{F}, \mathbb{P})$-martingale, we have $$\int_t^s e^{-\beta u}[r(u, X_u^{\pi^\prime}, \alpha_u, a_u^{\pi^\prime})-\hat{q}(u, X_u^{\pi^\prime}, \alpha_u, a_u^{\pi^\prime})+\hat{r}(u, X_u^{\pi^\prime}, \alpha_u, a_u^{\pi^\prime})] \mathrm{d} u$$ is an $(\mathbb{F}, \mathbb{P})$-martingale. Therefore, it is identically zero, which gives
$$r(u, X_u^{\pi^\prime}, \alpha_u, a_u^{\pi^\prime})-\hat{q}(u, X_u^{\pi^\prime}, \alpha_u, a_u^{\pi^\prime})+\hat{r}(u, X_u^{\pi^\prime}, \alpha_u, a_u^{\pi^\prime})=0,$$
for $\mathrm{d} u \otimes \mathrm{d} \mathbb{P}$-almost every point in $[t, T] \times \Omega$. By the continuous of the above expression, the full-support property of $\pi^{\prime}$, and the fact that the martingale condition holds for every initial condition $(t, x, i)$, we obtain $r(t, x, i, a)-\hat{q}(t, x, i, a)+\hat{r}(t, x, i, a)=0$ for every $(t, x, i, a)$. Hence,
\begin{align}\label{eq:q-r-identity}
\hat{q}(t, x, i, a)
&=r(t, x, i, a)+\hat{r}(t, x, i, a)
=H(t, x, i, a, \hat{J}(t, x, i))+\sum_{j=1}^L q_{i j} \hat{J}(t, x, j).
\end{align}
Substituting \eqref{eq:q-r-identity} into the constraint in \eqref{eq:The3.3-1} gives
\begin{align*}
\int_{\mathcal{A}}[H(t, x, i, a, \hat{J}(t, x, i))+\sum_{j=1}^L q_{i j} \hat{J}(t, x, j)+\gamma l_p( \pi(a \mid t, x, i))] \pi(a \mid t, x, i) \mathrm{d} a=0.
\end{align*}
We now identify $\hat{J}$ directly. Applying Itô's formula to $e^{-\beta s} \hat{J}(s, X_s^\pi, \alpha_s)$ along the state process under $\pi$, integrating from $t$ to $T$, taking expectations, and using the above identity, we obtain
\begin{align*}
\hat{J}(t, x, i)=\mathbb{E} \bigg[e^{-\beta(T-t)} \hat{J}(T, X_T^\pi, \alpha_T)+\int_t^T e^{-\beta(s-t)} \int_{\mathcal{A}}\big[r(s, X_s^\pi, \alpha_s, a)+\gamma l_p(\pi(a \mid s, X_s^\pi, \alpha_s))\big] \pi(a \mid s, X_s^\pi, \alpha_s) \mathrm{d} a \mathrm{d} s\bigg].
\end{align*}
The right-hand side is exactly $J(t, x, i ; \pi)$. Hence, $\hat{J}(t,x,i) \equiv J(t, x, i; \pi)$. Substituting this identity into the relation for $\hat{q}$ yields $\hat{q}(t,x,i,a) \equiv q(t,x,i,a; \pi)$.

Finally, if 
$\pi(a \mid t, x, i)=\big(\frac{p-1}{p \gamma}\big)^{\frac{1}{p-1}}(\hat{q}(t, x, i, a)+\frac{\gamma}{p-1}+\psi_p(t, x,i))_{+}^{\frac{1}{p-1}}, p \geq 1,$
satisfies $\int_{\mathcal{A}}\pi(a|t,x,i) \mathrm{d} a=1$, 
then $\pi=\mathcal{I}_p (\pi)$ where $\mathcal{I}_p(\cdot)$ is the map defined in \cref{thm:policy-improve}. This in turn implies $\pi(a \mid t, x, i)$ 
is an optimal policy and $\hat{J}(t,x,i)$ is the corresponding optimal value function.
\end{proof}

\section{q-Learning Algorithms under Tsallis Entropy}
In this section, we present learning algorithms derived from the martingale characterization of the q-function established in the previous section. We distinguish between two cases, based on whether the normalizing function can be explicitly computed and integrated.

\subsection{q-Learning algorithm when the normalizing function is available}\label{subsection:4.1}
In this subsection, we design q-learning algorithms to simultaneously learn and update the parameterized value function and the policy, utilizing the martingale condition established in Theorem \ref{thm:martingale-characterization}.

We first consider the case when the normalizing function is known or computable. Given a policy $\pi \in \Pi$, we parameterize the value function by a family of functions $J^\theta(\cdot, \cdot, \cdot)$, where $\theta \in \Theta \subset \mathbb{R}^{\mathsf{L}_\theta}$ and $\mathsf{L}_\theta$ is the parameter dimension. Similarly, we parameterize the q-function by $q^\zeta(\cdot, \cdot, \cdot, \cdot)$, where $\zeta \in \Psi \subset \mathbb{R}^{\mathsf{L}_\zeta}$ and $\mathsf{L}_\zeta$ is the parameter dimension. Consequently, the normalizing function $\psi^\zeta(t, x,i)$ is derived from the constraint
\begin{align}\label{eq:normalizing-fun}
\int_{\mathcal{A}}\left(\frac{p-1}{p \gamma}\right)^{\frac{1}{p-1}}\left(q^\zeta(t, x,i, a)+\psi^\zeta(t, x,i)\right)_{+}^{\frac{1}{p-1}} \mathrm{d} a=1.
\end{align}

Moreover, the approximators $J^\theta$ and $q^\zeta$ should also satisfy
\begin{align}\label{eq:approximators}
J^\theta(T, x, i)=h(x), \quad \int_{\mathcal{A}}[q^\zeta(t, x, i, a)+\gamma l_p(\pi^\zeta(a \mid t, x, i))] \pi^\zeta(a \mid t, x, i) \mathrm{d} a=0,
\end{align}
where the policy $\pi^\zeta$ is given by, for all $(t, x, i, a) \in[0, T] \times \mathbb{R}^d \times\mathcal{M} \times \mathcal{A}$,
$$
\pi^\zeta(a \mid t, x, i)=\left(\frac{p-1}{p \gamma}\right)^{\frac{1}{p-1}}\left(q^\zeta(t, x, i, a)+\psi^\zeta(t, x,i)\right)_{+}^{\frac{1}{p-1}}.
$$

Then, the learning task is to find the “optimal” (in some sense) parameters $\theta$ and $\zeta$. The key step in the algorithm design is to enforce the martingale condition stipulated in \cref{thm:policy-improve}.
By using the martingale orthogonality condition, it is enough to explore the solution $(\theta^*, \zeta^*)$ of the following martingale orthogonality equation system
\begin{align*}
&\Ex\left[\int_0^T \varrho_t\left(\mathrm{d} J^\theta(t, X_t^{\pi^\zeta}, \alpha_t)+r(t, X_t^{\pi^\zeta}, \alpha_t, a_t^{\pi^\zeta}) \mathrm{d} t-q^\zeta(t, X_t^{\pi^\zeta}, \alpha_t, a_t^{\pi^\zeta}) \mathrm{d} t-\beta J^\theta(t, X_t^{\pi^\zeta}, \alpha_t)\mathrm{d} t\right)\right]=0,\\
&\Ex\left[\int_0^T \varsigma_t\left(\mathrm{d} J^\theta(t, X_t^{\pi^\zeta},  \alpha_t)+r(t, X_t^{\pi^\zeta}, \alpha_t, a_t^{\pi^\zeta}) \mathrm{d} t-q^\zeta(t, X_t^{\pi^\zeta}, \alpha_t, a_t^{\pi^\zeta}) \mathrm{d} t-\beta J^\theta(t, X_t^{\pi^\zeta}, \alpha_t)\mathrm{d} t\right)\right]=0,
\end{align*}
where the test functions $\varrho=(\varrho_t)_{t \in[0, T]}, \varsigma=(\varsigma_t)_{t \in[0, T]}$ are $\mathbb{F}$-adapted stochastic processes. This can be implemented offline by using stochastic approximation to update parameters as
\begin{align}\label{eq:available}
\left\{\begin{array}{l}
\theta \leftarrow \theta+\eta_\theta \int_0^T \varrho_t\bigg(\mathrm{d} J^\theta(t, X_t^{\pi^\zeta},  \alpha_t)+r(t, X_t^{\pi^\zeta}, \alpha_t, a_t^{\pi^\zeta}) \mathrm{d} t-q^\zeta(t, X_t^{\pi^\zeta}, \alpha_t, a_t^{\pi^\zeta}) \mathrm{d} t-\beta J^\theta(t, X_t^{\pi^\zeta}, \alpha_t)\mathrm{d} t\bigg), \\
\zeta \leftarrow \zeta+\eta_\zeta \int_0^T \varsigma_t\bigg(\mathrm{d} J^\theta(t, X_t^{\pi^\zeta},  \alpha_t)+r(t, X_t^{\pi^\zeta}, \alpha_t, a_t^{\pi^\zeta}) \mathrm{d} t-q^\zeta(t, X_t^{\pi^\zeta}, \alpha_t, a_t^{\pi^\zeta}) \mathrm{d} t-\beta J^\theta(t, X_t^{\pi^\zeta}, \alpha_t)\mathrm{d} t\bigg),
\end{array}\right.
\end{align}
where $\eta_\theta$ and $\eta_\zeta$ are learning rates. In this paper, we choose the test functions in the conventional sense by
$$
\varrho_t=\frac{\partial J^\theta}{\partial \theta}(t, X_t^{\pi^\zeta}, \alpha_t), \quad \varsigma_t=\frac{\partial q^\zeta}{\partial \zeta}(t, X_t^{\pi^\zeta}, \alpha_t, a_t^{\pi^\zeta}).
$$
Based on the above updating rules, we present the pseudo-code of the offline q-learning algorithm in Algorithm \ref{Alg:Tsallis-q-Learning}.

\begin{algorithm}[htbp]
	\caption{\textbf{Offline q-Learning Algorithm with regime-switching}}
	\label{Alg:Tsallis-q-Learning}
	\hspace*{0.02in} {\bf Input:}
		Initial state-regime pair $(x, i)$, horizon $T$, time step $\Delta t$, number of episodes $N$, number of mesh grids $K$, the Markov Chain generator $Q$ and the state space $\mathcal{M}$, initial learning rates $\eta_\theta(\cdot), \eta_\zeta(\cdot)$(functions of the number of episodes), functional forms of parameterized value function $J^\theta(\cdot,\cdot,\cdot)$, and q-function $q^\zeta(\cdot,\cdot,\cdot,\cdot)$ satisfying \eqref{eq:approximators}, and temperature parameter $\gamma$.\\
		\hspace*{0.02in} {\bf Required Program:} environment simulator $(x^{\prime}, r, j)=$ Environment $_{\Delta t}(t, x, i, a, Q)$ that takes current time-state-regime-action quadruple $(t, x, i, a)$ and the $Q$ matrix as inputs, and generates next state $x^{\prime}$, instantaneous reward $r$, and next regime $j$ (driven by the Markov Chain $Q$) at time $t+\Delta t$ as outputs. \\
		\hspace*{0.02in} {\bf Learning Procedure:}

        \begin{algorithmic}[1]
            \State Initialize $\theta, \zeta$, and $n=1$.
			\While{$n<N$}
			      \State Initialize $\jmath = 0$.  Observe initial state $x$, regime $i$, and store $x_{t_0}\leftarrow x, \alpha_{t_0}\leftarrow i$.
				\While{$\jmath < K$}
				    \State Generate action $a_{t_\jmath} \sim \pi^\zeta(\cdot \mid t_\jmath, x_{t_\jmath} \alpha_{t_\jmath})$.
                    \State Apply $a_{t_\jmath} $ to environment simulator $(x, r, j)=$ Environment  $_{\Delta t}(t_\jmath, x_{t_\jmath}, \alpha_{t_\jmath}, a_{t_\jmath}, Q)$.
                    \State Store $x_{t_{\jmath+1}} \leftarrow x$, $\alpha_{t_{\jmath+1}} \leftarrow j$, and $r_{t_{\jmath+1}} \leftarrow  r$.
                    \State Update time index: $t_{\jmath+1} \leftarrow t_\jmath + \Delta t$, and step counter: $\jmath \leftarrow \jmath+1$.
                \EndWhile{}
                \State For every $k=0,1,...,K-1$, compute
                    \begin{align*}                              G_{k}=&J^\theta\left(t_{k+1},x_{t_{k+1}},\alpha_{t_{k+1}}\right)-J^\theta\left(t_{k},x_{t_k},\alpha_{t_{k}}\right)+r_{t_k}\Delta t-q^{\zeta}\left(t_{k},x_{t_{k}},\alpha_{t_{k}}, a_{t_{k}}\right) \Delta t-\beta J^\theta\left(t_{k+1},x_{t_{k+1}},\alpha_{t_{k+1}}\right)\Delta t.
                    \end{align*}
				\State Update $\theta$ and $\zeta$ by
                    \begin{align*}
                        \theta \leftarrow \theta+\eta_\theta\sum_{k=0}^{K-1} \frac{\partial J^\theta}{\partial \theta}\left(t_k,x_{t_k},\alpha_{t_{k}}\right)G_{k}, \quad
                        \zeta \leftarrow \zeta+\eta_\zeta\sum_{k=0}^{K-1}   \frac{\partial q^\zeta}{\partial \zeta}\left(t_k,x_{t_k},\alpha_{t_{k}},a_{t_k}\right)G_{k}.
                    \end{align*}
               \State Update $n \leftarrow n+1$.
              \EndWhile{}
    \end{algorithmic}
\end{algorithm}

\subsection{q-Learning algorithm when the normalizing function is unavailable}

In this subsection, we address the case when the normalizing function does not admit an explicit form. 
We still parameterize the value function by $J^\theta(\cdot, \cdot, \cdot)$, where $\theta \in \Theta \subset \mathbb{R}^{\mathsf{L}_\theta}$ (with dimension $\mathsf{L}_{\theta}$), and the q-function by $q^\zeta(\cdot, \cdot, \cdot, \cdot)$, where $\zeta \in \Psi \subset \mathbb{R}^{\mathsf{L}_\zeta}$ (with dimension $\mathsf{L}_\zeta$). However, since $\psi^\zeta(t, x, i)$ cannot be derived from \eqref{eq:normalizing-fun}, we instead parameterize the policy by a family of functions $\pi^\chi(\cdot)$, where $\chi \in \Upsilon \subset \mathbb{R}^{\mathsf{L}_\chi}$ (with dimension $\mathsf{L}_\chi$). Furthermore, the approximators $J^\theta$ and $\pi^\chi$ must satisfy the terminal condition $J^\theta(T, x, i)=h(x)$. We then define the function $F:[0, T] \times \mathbb{R}^d\times \mathcal{M} \times \mathcal{P}(\mathcal{A}) \times \mathcal{P}(\mathcal{A}) \mapsto \mathbb{R}$ as
\begin{align*}
F(t, x, i ; \pi^{\prime}, \pi):=\int_{\mathcal{A}}[q(t, x, i, a ; \pi)+\gamma l_p(\pi^{\prime}(a \mid t, x, i))] \pi^{\prime}(a \mid t, x, i) \mathrm{d} a .
\end{align*}
Then, we can devise an Actor-Critic q-learning algorithm to learn the q-function and the optimal policy in an alternating manner. For the Actor-step (or policy improvement step), we update the policy $\pi^\chi$ by maximizing the function $F(t, x, i; \pi^{\chi^{\prime}}, \pi^\chi)$ that
$$
\max_{\chi^{\prime} \in \Upsilon} F(t, x, i ; \pi^{\chi^{\prime}}, \pi^\chi)=\max_{\chi^{\prime} \in \Upsilon} \int_{\mathcal{A}}[q(t, x, i, a; \pi^\chi)+\gamma l_p(\pi^{\chi^{\prime}}(a \mid t, x, i))] \pi^{\chi^{\prime}}(a \mid t, x, i) \mathrm{d} a.
$$

In fact, we have the following result, which is a direct consequence of Theorem \ref{thm:policy-improve}. \citet{bo2024continuus} proved such a theorem while studying continuous-time q-learning in jump diffusion models. We then apply this method to the Markov regime-switching system.
\begin{lemma} \label{lem:J}
    Given $(t, x, i) \in[0, T] \times \mathbb{R}^d \times \mathcal{M}$ and $\pi, \pi^{\prime} \in \Pi$, if it holds that $F(t, x, i; \pi^{\prime}, \pi) \geq F(t, x, i; \pi, \pi)$, then $J(t, x, i; \pi^{\prime}) \geq J(t, x, i; \pi)$.
\end{lemma}

Moreover, while the q-learning method based on Theorem \ref{thm:martingale-characterization} requires the policy function $\pi^\chi$ to satisfy $\pi^\chi \in \mathcal{P}(\mathcal{A})$ and the consistency condition \eqref{eq:The3.3-1}, we relax these constraints here. Instead, we consider the following maximization problem, for $w_1, w_2 \geq 0$,
$$
\max_{\chi^{\prime} \in \Upsilon}\bigg[F(t, x, i ; \pi^{\chi^{\prime}}, \pi^\chi)-w_1 F^2(t, x, i; \pi^{\chi^{\prime}}, \pi^{\chi^{\prime}})-w_2\left(\int_{\mathcal{A}} \pi^{\chi^{\prime}}(a\mid t, x, i) \mathrm{d} a-1\right)^2\bigg].
$$
By a direct calculation, we obtain
\begin{align*}
\frac{\partial F(t, x, i; \pi^{\chi^{\prime}}, \pi^\chi)}{\partial \chi^{\prime}}
=&\int_{\mathcal{A}}\left(q(t, x, i, a ; \pi^\chi)+\gamma l_p(\pi^{\chi^{\prime}}(a \mid  t, x, i))\right) \frac{\partial \pi^{\chi^{\prime}}(a \mid t, x, i)}{\partial \chi^{\prime}} \mathrm{d} a\\
&+\gamma\int_{\mathcal{A}} l_p^{\prime}(\pi^{\chi^{\prime}}(a \mid  t, x, i)) \frac{\partial \pi^{\chi^{\prime}}(a \mid t, x, i)}{\partial \chi^{\prime}} \pi^{\chi^{\prime}}(a \mid t, x, i) \mathrm{d} a \\
=&\int_{\mathcal{A}}\left(q(t, x, i, a; \pi^\chi)+\gamma l_p(\pi^{\chi^{\prime}}(a \mid t, x, i))\right) \frac{\partial \ln \pi^{\chi^{\prime}}(a \mid t, x, i)}{\partial \chi^{\prime}} \pi^{\chi^{\prime}}(a \mid t, x, i) \mathrm{d} a \\
&+\gamma \int_{\mathcal{A}} l_p^{\prime}(\pi^{\chi^{\prime}}(a \mid t, x, i)) \frac{\partial \pi^{\chi^{\prime}}(a \mid t, x, i)}{\partial \chi} \pi^{\chi^{\prime}}(a \mid t, x, i) \mathrm{d} a.
\end{align*}
Hence, we can update $\chi$ by using stochastic gradient descent:
\begin{align}\label{eq:update-rule-F}
\chi \leftarrow \chi & +\eta_\chi\Bigg\{\int_{0}^{T} \bigg[\left(q(t, X_t, \alpha_t, a^{\pi^\chi} ; \pi^\chi)+\gamma l_p(\pi^{\chi}(a^{\pi^\chi} \mid  t, X_t, \alpha_t))\right) \frac{\partial \ln \pi^\chi(a^{\pi^\chi} \mid  t, X_t, \alpha_t)}{\partial \chi} \notag\\
&+\gamma l_p^{\prime}(\pi^\chi(a^{\pi^\chi} \mid  t, X_t, \alpha_t)) \frac{\partial \pi^\chi(a^{\pi^\chi} \mid  t, X_t, \alpha_t)}{\partial \chi}\bigg] \mathrm{d} t-2 w_1 \int_0^T F(t, X_t, \alpha_t ; \pi^\chi, \pi^\chi) \frac{\partial F(t, X_t, \alpha_t ; \pi^\chi, \pi^\chi)}{\partial \chi} \mathrm{d} t \notag\\
&-2 w_2 \int_0^T\left(\int_{\mathcal{A}} \pi^\chi(a \mid  t, X_t, \alpha_t) \mathrm{d} a-1\right) \int_{\mathcal{A}} \frac{\partial}{\partial \chi}\pi^\chi(a \mid t, X_t, \alpha_t) \mathrm{d} a \mathrm{d} t\bigg\},
\end{align}
where $\eta_\chi$ is the learning rate.

Next, for the Critic-step (or the policy evaluation step), we adopt the same parameter updating rules for the value function and q-function as described in \eqref{eq:available} in the previous algorithm (Subsection \ref{subsection:4.1}). We present the pseudo-code
in Algorithm \ref{Alg:Tsallis-q-Learning-normalizing-unavailable-F}.

\begin{algorithm}[htbp]	\caption{\textbf{Offline q-Learning Algorithm with regime-switching (Normalizing Function Unavailable)}}
		    \label{Alg:Tsallis-q-Learning-normalizing-unavailable-F}
			\hspace*{0.02in} {\bf Input:}
			Initial state-regime pair $(x, i)$, horizon $T$, time step $\Delta t$, number of episodes $N$, number of mesh grids $K$, the Markov Chain generator $Q$ and the state space $\mathcal{M}$, initial learning rates $\eta_\theta(\cdot), \eta_\zeta(\cdot), \eta_\chi(\cdot)$(functions of the number of episodes), functional forms of parameterized value function $J^\theta(\cdot,\cdot,\cdot)$, and q-function $q^\zeta(\cdot,\cdot,\cdot,\cdot)$ satisfying \eqref{eq:approximators}, and temperature parameter $\gamma$.\\
			\hspace*{0.02in} {\bf Required Program:} environment simulator $(x^{\prime}, r, j)=$ Environment $_{\Delta t}(t, x, i, a, Q)$ that takes current time-state-regime-action quadruple $(t, x, i, a)$ and the $Q$ matrix as inputs, and generates next state $x^{\prime}$, instantaneous reward $r$, and next regime $j$ (driven by the Markov Chain $Q$) at time $t+\Delta t$ as outputs. \\
			\hspace*{0.02in} {\bf Learning Procedure:}
			\begin{algorithmic}[1]
\State Initialize $\theta,\zeta,\chi$, and $n=1$.
				\While{$n<N$}
				\State Initialize $\jmath = 0$. Observe initial state $x$, regime $i$, and store $x_{t_0}\leftarrow x, \alpha_{t_0}\leftarrow i$.
				\While{$\jmath < K$}
				    \State Generate action $a_{t_\jmath} \sim \pi^\chi(\cdot \mid t_\jmath, x_{t_\jmath} \alpha_{t_\jmath})$. Apply $a_{t_\jmath} $ to environment simulator $(x, r, j)=$ \State Environment  $_{\Delta t}(t_\jmath, x_{t_\jmath}, \alpha_{t_\jmath}, a_{t_\jmath}, Q)$. Store $x_{t_{\jmath+1}} \leftarrow x$, $\alpha_{t_{\jmath+1}} \leftarrow j$, and $ r_{t_{\jmath+1}} \leftarrow  r$. \State Update time index: $t_{\jmath+1} \leftarrow t_\jmath + \Delta t$, and step counter: $\jmath \leftarrow \jmath+1$.
                \EndWhile{}
                \State  For every $k=0,1,...,K-1$, compute
\begin{align*}
 G_{k}=J^\theta\left(t_{k+1},x_{t_{k+1}},\alpha_{t_{k+1}}\right)-J^\theta\left(t_{k},x_{t_k},\alpha_{t_{k}}\right)+r_{t_k}\Delta t-q^{\zeta}\left(t_{k},x_{t_{k}},\alpha_{t_{k}},a_{t_{k}}\right) \Delta t
 -\beta J^\theta\left(t_{k+1},x_{t_{k+1}},\alpha_{t_{k+1}}\right)\Delta t.
\end{align*}
				\State For the Critic (policy evaluation) step, update $\theta$ and $\zeta$ by
\begin{align*}
\theta &\leftarrow \theta+\eta_\theta \sum_{k=0}^{K-1} \frac{\partial J^\theta}{\partial \theta}\left(t_k,x_{t_k},\alpha_{t_k}\right)G_{k}, \quad
\zeta \leftarrow \zeta+\eta_\zeta \sum_{k=0}^{K-1}   \frac{\partial q^\zeta}{\partial \zeta}\left(t_k,x_{t_k}, \alpha_{t_k},a_{t_k}\right)G_{k}.
\end{align*}
\State For the Actor (policy improvement) step, update $\chi$ (using the updated $\theta$ and $\zeta$) by
\begin{align*}
\chi\leftarrow\chi&+\eta_\chi\Bigg\{\sum_{k=0}^{K-1}\bigg[(q^{\zeta}(t_k,x_{t_k},\alpha_{t_k},a_{t_k})+\gamma l(\pi^{\chi}(a_{t_k}\mid t_k,x_{t_k},\alpha_{t_k})) \frac{\partial \ln \pi^{\chi} (a_{t_k}\mid t_k,x_{t_k},\alpha_{t_k})}{\partial \chi}\\
&+\gamma l'(\pi^{\chi}(a_{t_k}\mid t_k,x_{t_k},\alpha_{t_k})) \frac{\partial \pi^{\chi} (a_{t_k}\mid t_k,x_{t_k},\alpha_{t_k})}{\partial \chi}\bigg]-2w_1 \sum_{k=0}^{K-1}F(t_k,x_{t_k},\alpha_{t_k};\pi^{\chi},\pi^{\chi})\frac{\partial F(t_k,x_{t_k},\alpha_{t_k};\pi^{\chi},\pi^{\chi})}{\partial \chi}\\
&-2w_2\sum_{k=0}^{K-1} \left(\int_{\mathcal{A}} \pi^{\chi}(a_{t_k}\mid t_k,x_{t_k},\alpha_{t_k})\mathrm{d}a-1\right)\int_{\mathcal{A}} \frac{\partial \pi^{\chi}}{\partial \chi}(a_{t_k}\mid t_k,x_{t_k},\alpha_{t_k})\mathrm{d}a\Bigg\}.
\end{align*}
               \State   Update $n \leftarrow n+1$.
              \EndWhile{}
    \end{algorithmic}
\end{algorithm}


To overcome the difficulty of an unavailable normalizing constant in soft Q-learning, \citet{Haarnoja2018Soft} also introduced a general method of using a family of stochastic policies whose densities can be easily computed to approximate $\pi^{\prime}$. Specifically, denote by $\{\boldsymbol{\pi}^\phi(\cdot \mid t, x, i)\}_{\phi \in \Phi}$ the family of density functions of some tractable distributions such as Gaussians. The learning procedure starts with a policy $\pi^\phi$ from this family. When $p=1$, the objective is to project the target policy $\frac{\exp{\frac{1}{\gamma}}q(t, x, i, \cdot; \pi^\phi)}{\int_{\mathcal{A}}\exp{\frac{1}{\gamma}}q(t, x, i, a;\pi^\phi)\mathrm{d} a}$ by minimizing
\begin{align*}
&\min_{\phi^{\prime} \in \Phi} D_{K L}\left(\pi^{\phi^{\prime}}(\cdot \mid t, x, i) \| \frac{\exp{\frac{1}{\gamma}}q(t, x, i, \cdot; \pi^\phi)}{\int_{\mathcal{A}}\exp{\frac{1}{\gamma}}q(t, x, i, a;\pi^\phi)\mathrm{d} a}\right) =\min_{\phi^{\prime} \in \Phi} D_{K L}\left(\pi^{\phi^{\prime}}(\cdot \mid t, x, i) \| \exp{\frac{1}{\gamma}}q(t, x, i, \cdot; \pi^\phi)\right),
\end{align*}
where $D_{K L}(f \| g):=\int_{\mathcal{A}} \log \frac{f(a)}{g(a)} f(a) \mathrm{d} a$ is the Kullback-Leibler (KL) divergence of two positive functions $f, g$ with the same support on $\mathcal{A}$, where $f \in \mathcal{P}(\mathcal{A})$ is a probability density function on $\mathcal{A}$.

We present a policy improvement theorem, similar to Theorem 10 in \citet{jia2023q}, that uses the KL divergence to optimize the policy without explicitly computing the normalization constant.
\begin{theorem} \label{thm:KL}
Given $(t, x, i)\in[0, T]\times\mathbb{R}^d\times\mathcal{M}$, if two policies $\pi \in \Pi$ and $\pi^{\prime} \in \Pi$ satisfy
\begin{align*}
& D_{K L}\left(\pi^\prime(\cdot \mid t, x, i)\| \exp \{\frac{1}{\gamma} H(t, x, i, \cdot, J(t, x, i ; \pi))\}\right)
\leq D_{K L}\left(\pi(\cdot \mid t, x, i) \| \exp \{\frac{1}{\gamma} H(t, x, i, \cdot, J(t, x, i ; \pi))\}\right),
\end{align*}
then $J(t, x, i; \pi^\prime) \geq J(t, x, i ; \pi)$.
\end{theorem}

Theorem \ref{thm:KL} presents a general result for comparing any two given policies, regardless of whether they belong to a tractable family of densities.
Therefore, similar to the discussion in Theorem 10 of \citet{jia2023q}, we can update $\phi$ incrementally at each step
\begin{align}\label{eq:KL-update}
\phi\leftarrow \phi - \eta_\phi \left[\log \pi^{\phi}(a^{\pi^{\phi}} \mid t, X_t, \alpha_t)-\frac{1}{\gamma} q(t, X_t, \alpha_t, a^{\pi^{\phi}}; \pi^\phi)\right]\frac{\partial}{\partial \phi}\log \pi^{\phi}(a^{\pi^{\phi}} \mid t, X_t, \alpha_t),
\end{align}
where $\eta_\phi$ is the learning rate. 
This approach is known as variational inference (or optimal Gaussian approximation) (see \citet{Murphy2012Machine}). Essentially, it simplifies analysis by approximating complex distributions, functions, or problems with a Gaussian (normal) distribution. However, this approach is not applicable to all policy updates when the normalizing function is unavailable. More specifically, this update method is suitable only when the optimal policy distribution is unimodal. Therefore, we provide this method only as a solution for cases in which the optimal policy distribution is unimodal.

Analysis shows that for the order of Tsallis entropy $p>1$, the objective is to project the target policy
\begin{align*}
\pi^\phi(a \mid t, x, i)=\left(\frac{p-1}{p \gamma}\right)^{\frac{1}{p-1}}(q(t, x, i, \cdot;\pi^\phi)+\psi(t, x, i;\pi^\phi))_{+}^{\frac{1}{p-1}},
\end{align*}
by minimizing $\min_{\phi^{\prime} \in \Phi} D_{K L}(\pi^{\phi^{\prime}}(\cdot \mid t, x, i) \| \pi^\phi(a \mid t, x, i),$
where $\psi(t, x, i;\pi^\phi)$ satisfies $\int_{\mathcal{A}}\pi^\phi(a \mid t, x, i)\mathrm{d} a=1$.
It is worth noting that this is the forward KL divergence (see \citet{Murphy2012Machine}), which requires the target policy $\pi^\phi(a \mid t, x, i)$ to be unimodal or strictly concave. For $p=1$, the target policy $\frac{\exp{\{\frac{1}{\gamma}}q(t, x, i, \cdot; \pi^\phi)\}}{\int_{\mathcal{A}}\exp{\{\frac{1}{\gamma}}q(t, x, i, a;\pi^\phi)\}\mathrm{d} a}$ is a Gaussian policy, which is clearly unimodal. To ensure the uniqueness of the optimal policy distribution in most standard optimal control models, the target policy $\pi^\phi(a \mid t, x, i)$ is generally required to be strictly concave with respect to $a$. Therefore, we exclusively consider cases where $\pi^\phi(a \mid t, x, i)$ is strictly concave in $a$.

Define the multivariate normal distribution $\boldsymbol{\pi}^{\phi^{\prime}}(\cdot \mid t, x, i)=\mathcal{N}(\mu^{\phi^{\prime}}(t, x, i), \Sigma^{\phi^{\prime}}(t, x, i))$, where $\mu^{\phi^{\prime}}(t, x, i) \in \mathbb{R}^m$ and $\Sigma^{\phi^{\prime}}(t, x, i) \in \mathbb{S}_{++}^m$.
We now consider utilizing the optimal Gaussian approximation (OGA) to solve the policy update problem $\min_{\phi^{\prime} \in \Phi} D_{K L}(\pi^{\phi^\prime}(a \mid t, x, i) \| \pi^\phi(a \mid t, x, i))$. The core method employs a second-order Taylor expansion to rewrite the complex KL-divergence minimization problem as a local quadratic program with a closed-form solution. As the minimization of $D_{K L}(\pi^{\phi^\prime}(a \mid t, x, i) \parallel \pi^\phi(a \mid t, x, i))$ is generally intractable analytically, OGA yields an approximate closed-form solution via the second-order Taylor expansion of the target log density $\log \pi^\phi(a \mid t, x, i)$ near the current mean $\mu$. Minimizing $D_{K L}(\pi^{\phi^\prime}(a \mid t, x, i) \| \pi^\phi(a \mid t, x, i))$ is equivalent to finding a Gaussian distribution $\pi^{\phi^\prime}(a \mid t, x, i)$ such that its log density $\log \pi^{\phi^\prime}(a \mid t, x, i)$ matches the second order Taylor expansion of the target log density $\log \pi^\phi(a \mid t, x, i)$ at the current mean $\mu$.

\section{Applications and Numerical Examples}
Consider an investor who manages a portfolio with an investment horizon $T>0$. For clarity of presentation, the market consists of one risky asset, the stock $\{S_t\}_{t \in[0, T]}$, and one risk-free asset, the bond $\{B_t\}_{t \in[0, T]}$. Let $\{W_t\}_{t \in[0, T]}$ represent a one-dimensional Brownian Motion defined on the filtered probability space $(\Omega, \mathcal{F}, \mathbb{F}, \mathbb{P})$ that satisfies the usual conditions. For any time $t \in[0, T], \alpha_t$ takes a value from the set $\mathcal{M}$.
The dynamics of the stock and the bond are driven by two stochastic processes
\begin{align*}
& \mathrm{d} S_t=S_t\{\mu(t, \alpha_t) \mathrm{d} t+\sigma(t, \alpha_t) \mathrm{d} W_t\}, \text { with } S_0>0, \\
& \mathrm{d} B_t=r(t, \alpha_t) B_t \mathrm{d} t, \quad\text{with} \quad B_0>0,
\end{align*}
where $\mu(t, \alpha_t) \in \mathbb{R}$ and $\sigma(t, \alpha_t) >0$ are the mean and volatility of the stock return, respectively. And $r(t, \alpha_t) \in \mathbb{R}_{+}$ is the risk-free interest rate at time $t$ in market regime $\alpha_t$. Here, $r(t,\alpha_t)$ denotes the risk-free interest rate and should be distinguished from the running reward function $r(t,x,i,a)$ introduced in the general framework above. 
At each time $t \in [0, T]$, let $X_t$ denote the value of the investor's portfolio. The investor reallocates their portfolio by investing $a_t$ in stocks and $X_t - a_t$ in bonds. Under the self-financing assumption, the portfolio value process can be derived as
\begin{align}\label{eq:MV-portfolio-value}
\mathrm{d} X_t^a
& =[r(t, \alpha_t) X_t^a+\rho(t, \alpha_t) \sigma(t, \alpha_t) a_t] \mathrm{d} t+\sigma(t, \alpha_t) a_t \mathrm{d} W_t,
\end{align}
given the initial portfolio value $X_0=x_0>0$ and the initial regime $\alpha_0=i$. Here, $\rho(t, \alpha_t):=\sigma^{-1}(t, \alpha_t)(\mu(t, \alpha_t)-r(t, \alpha_t))$ is the Sharpe ratio, and $\{X_t^a\}_{t \in[0, T]}$ with a superscript $a$ represents the portfolio value process that follows the control policy $a:=\{a_t\}_{t \in[0, T]}$.
The classical continuous-time MV model aims to solve the following constrained optimization problem
\begin{equation}\label{eq:classical-MV-problem}
 \min_a \operatorname{Var}[X_T^a],\quad \text { subject to } \Ex[X_T^a]=z,
\end{equation}
where $\{X_t^a, 0 \leq t \leq T\}$ satisfies the dynamics \eqref{eq:MV-portfolio-value} under the investment strategy (portfolio) $a$, and $z \in \mathbb{R}$ is an investment target set at $t=0$ as the desired mean payoff at the end of the investment horizon $[0, T]$. To solve \eqref{eq:classical-MV-problem}, one first transforms it into an unconstrained problem by applying a Lagrange multiplier $w$
\begin{equation*}
\min_a \{\Ex[(X_T^a)^2]-z^2-2 w(\Ex[X_T^a]-z)\}=\min_a \{\Ex[(X_T^a-w)^2]-(w-z)^2\}.
\end{equation*}
Equivalently, for each fixed $w$, the same optimal investment strategy can be obtained from the following maximization problem:
\begin{equation*}
\max_a \{-\Ex[(X_T^a-w)^2] +(w-z)^2\}.
\end{equation*}
This problem can be solved analytically, and its $a^*=\{a_t^*, 0 \leq t \leq T\}$ depends on $w$. Then the original constraint $\Ex[X_T^{a^*}]=z$ determines the value of $w$. We refer the reader to \citet{zhou2000Continuous} for a detailed derivation.

We now present the entropy-regularized exploratory counterpart of this problem.
Under an admissible policy $\pi \in \Pi$ and $\mathcal{A}=\mathbb{R}$, the exploratory counterpart of the dynamics \eqref{eq:MV-portfolio-value} is governed by the following SDE
\begin{equation*}
    \begin{aligned}
        \mathrm{d} \tilde{X}_t^\pi=&\bigg[r(t, \alpha_t) \tilde{X}_t^\pi+\int_{\mathcal{A}} \rho(t, \alpha_t) \sigma(t, \alpha_t)a\cdot\pi(a \mid t, \tilde{X}_t^\pi, \alpha_t) \mathrm{d} a\bigg] \mathrm{d} t+\bigg(\sqrt{\int_{\mathcal{A}} \sigma^2(t, \alpha_t)a^2\cdot\pi(a \mid t, \tilde{X}_t^\pi, \alpha_t) \mathrm{d} a}\bigg) \mathrm{d} W_t.
   \end{aligned}
\end{equation*}
For any $(t, x, i) \in[0, T] \times \mathbb{R} \times \mathcal{M}$, we define the value function $J$ as
\begin{equation*}
    J^\pi(t, x, i;w):=\Ex\bigg[-(\tilde{X}_T^\pi-w)^2+\gamma \int_t^T \int_{\mathcal{A}}l_p(\pi(a \mid s, \tilde{X}_s^\pi,\alpha_s))\pi(a \mid s, \tilde{X}_s^\pi,\alpha_s) \mathrm{d} a \mathrm{d} s\mid \tilde{X}_t^\pi=x, \alpha_t=i\bigg]+(w-z)^2. \\
\end{equation*}
Accordingly, the optimal value function $V$ is given by
$V(t, x, i;w):=\sup_{\pi \in \Pi}J^\pi(t, x, i;w)$.
By applying Itô's lemma and the dynamic programming principle, we obtain the HJB equation
\begin{equation}\label{eq:MV-HJB}
    \begin{aligned}
        &\sup_{\pi \in \Pi}\int_{\mathcal{A}}[H(t, x, i, a, V)+ \sum_{j=1}^L q_{ij}V(t, x, j)+\gamma l_p(\pi(a \mid t, x, i))]\pi(a \mid t, x, i) \mathrm{d} a=0,
    \end{aligned}
\end{equation}
with the terminal condition $V(T, x, i;w)=-(x-w)^2+(w-z)^2$. Here,
\begin{align}\label{eq:MV-Hamiltonian}
   H(t, x, i, a, V) = V_t(t, x, i;w) + \big[r(t, i) x + \rho(t, i) \sigma(t, i) a\big] V_x(t, x, i;w) + \frac{1}{2} \sigma^2(t, i) a^2 V_{xx}(t, x, i;w),
\end{align}
where the running reward $r(t,x,i,a)$ in the general framework and the discount rate $\beta$ are both set equal to zero.
Following a similar argument to that in Section \ref{subsection:HJB}, we can solve the (constrained) optimization problem in the HJB equation \eqref{eq:MV-HJB} to obtain a feedback (distributional) control whose density function is given by
\begin{align}\label{eq:MV-optimal-policy}
\pi_p^*(a \mid t, x, i)=
\begin{cases}
\left(\frac{p-1}{p \gamma}\right)^{\frac{1}{p-1}}(H(t, x, i, a, V)+\sum_{j=1}^Lq_{ij}V(t, x, j;w)+\psi(t, x, i))_{+}^{\frac{1}{p-1}},& p>1,\\
\mathcal{N}\left(a \left\lvert\,-\frac{\rho(t, i)V_x(t, x, i;w)}{\sigma(t, i)V_{x x}(t, x, i;w)}\right., -\frac{\gamma}{\sigma(t, i)^2 V_{x x}(t, x, i;w)}\right),& p=1,
\end{cases}
\end{align}
where $V_{x x}(t, x, i;w)<0$ and the normalizing function $\psi(t, x, i)$ is determined by $\int_{\mathcal{A}}\pi_p^*(a \mid t, x, i)\mathrm{d} a=1$.

Motivated by the quadratic terminal condition and the linear-quadratic structure of the wealth dynamics, we seek the value function in the quadratic form
\begin{align}\label{eq:MV-optimal-value-func-m1m2m3}
V(t, x, i;w)=A(t, i)[x+w B(t, i)]^2+w^2 C(t, i)+D(t, i)+(w-z)^2, \quad t\in[0, T],
\end{align}
where $A(t, i)<0$, $A(T,i)=-1$, $B(T,i)=-1$, $C(T,i)=0$, $D(T,i)=0$. 
Consequently, for $p > 1$, the optimal policy in \eqref{eq:MV-optimal-policy} can be rewritten as
\begin{equation}\label{eq:MV-optimal-policy-m1m2m3-p2}
\begin{aligned}
\widehat{\pi}_p(a \mid t, x, i)
=\bigg(\frac{p-1}{p \gamma}\bigg)^{\frac{1}{p-1}}\bigg(2 A&(t, i)\rho(t, i) \sigma(t, i)[x+ w B(t, i)]  a + A(t, i) \sigma^2(t, i) a^2+\sum_{j=1}^L q_{i j}V(t, x, j;w)\\
&+ \psi(t, x, i)\bigg)_{+}^{\frac{1}{p-1}} ,
\end{aligned}
\end{equation}
where $\psi(t, x, i)$ is determined by the constraint $\int_{\mathcal{A}} \widehat{\pi}_p(a \mid t, x, i)\mathrm{d} a=1$. For the Shannon entropy case ($p = 1$), the optimal policy reduces to
\begin{align*}
\widehat{\pi}_1(a \mid t, x, i)
=\mathcal{N}\left(a \left\lvert\,-\frac{\rho(t, i)(x+wB(t, i))}{\sigma(t, i)}\right., -\frac{\gamma}{2\sigma(t, i)^2 A(t, i)}\right).
\end{align*}
Moreover, at the initial time $t = 0$, the optimal Lagrange multiplier is given by
\begin{align}\label{eq:MV-optimal-w}
w^*=\frac{z-A(0, i) B(0, i) x_0}{A(0, i)B(0, i)^2+C(0, i)+1}.
\end{align}
In accordance with Definition \ref{def:q-corollary}, we are now well-positioned to characterize the q-function:
\begin{align}\label{eq:MV-q-func}
q^\pi(t, x, i, a; w)
=&J^\pi_{t}(t, x, i; w)+J^\pi_{x}(t, x, i; w)\left[r(t, i) x+\rho(t, i) \sigma(t, i)a\right]\notag\\
&+\frac{1}{2}J^\pi_{x x}(t, x, i; w) \sigma^2(t, i)a^2 +\sum_{j=1}^L q_{ij}J^\pi(t, x, j; w).
\end{align}
To solve Eq.\eqref{eq:MV-HJB}, we investigate the cases $p=1$ and $p>1$ separately, since evaluating the policy's required moments is straightforward for the Gaussian case ($p=1$) but more difficult for $p>1$ due to the implicit normalization function. For simplicity, we consider a two-state framework ($\mathcal{M} = \{1, 2\}$) to characterize two distinct market regimes in the following numerical experiment.

\textbf{The case of $p = 1$. }
In this case, the first and second moments of the Gaussian policy $\widehat{\pi}_1(a \mid t, x, i)$ (defined in \eqref{eq:MV-optimal-policy}) are explicitly given by
\begin{align*}
\Ex_{\widehat{\pi}_1}[a] =-\frac{\rho(t, i) V_x(t, x, i)}{\sigma(t, i) V_{xx}(t, x, i)},\quad
\Ex_{\widehat{\pi}_1}[a^2] = \left(\frac{\rho(t, i) V_x(t, x, i)}{\sigma(t, i) V_{xx}(t, x, i)}\right)^2 - \frac{\gamma}{\sigma^2(t, i) V_{xx}(t, x, i)}.
\end{align*}
Substituting these moments into \eqref{eq:MV-HJB} and recalling \eqref{eq:MV-optimal-value-func-m1m2m3} yields
\begin{align}\label{eq:equs-of-m1m2m3}
\begin{cases}
A_t(t, i)=[\rho^2(t, i)-2r(t,i)]A(t, i)-\sum_{j=1}^2 q_{i j} A(t, j),& A(T, i) = -1,\\
B_t(t, i)=r(t, i)B(t, i)-\frac{1}{A(t, i)} \sum_{j=1}^2 q_{i j} A(t, j)(B(t, j)-B(t, i)), & B(T, i) = -1,\\
C_t(t, i)=-\sum_{j=1}^2 q_{ij}C(t, j)-\sum_{j=1}^2 q_{ij}A(t, j)(B(t, j)-B(t, i))^2, & C(T, i) = 0,\\
D_t(t, i)=-\sum_{j=1}^2 q_{ij}D(t, j)-\frac{\gamma}{2} \log \left(\frac{\pi \gamma}{-\sigma^2(t, i) A(t, i)}\right),& D(T, i) = 0.
\end{cases}
\end{align}
Consequently, the corresponding q-function defined in \eqref{eq:MV-q-func} can be expressed as
\begin{align*}
q^{\widehat{\pi}_1}(t, x, i, a; w)
=&[A_t(t, i) +A(t, i)(2r(t, i)- \rho^2(t, i))][x+w B(t, i)]^2\notag\\
&+ 2w A(t, i) [B_t(t, i)-r(t, i)B(t, i)] [x+w B(t, i)]+ w^2 C_t(t, i) + D_t(t, i)\notag\\
& + A(t, i)\sigma^2(t, i)\bigg[a+\frac{\rho(t, i)(x+wB(t, i))}{\sigma(t, i)}\bigg]^2 + \sum_{j=1}^2 q_{ij}J^{\widehat{\pi}_1}(t, x, j; w).
\end{align*}
To simplify the notation, we introduce the column vector $\mathbf{F}(t) = (F(t, 1), F(t, 2))^\top$ for any function $F \in \{A, B, C, D\}$, with its time derivative naturally denoted by $\mathbf{F}_t(t)$.
To facilitate mathematical tractability, we rewrite the differential system \eqref{eq:equs-of-m1m2m3} into a compact vector-matrix ODE format. Assuming time-homogeneous market parameters (i.e., $\rho(t, i) \equiv \rho_i$, $r(t, i) \equiv r_i$, and $\sigma(t, i) \equiv \sigma_i$), we obtain
\begin{align}\label{eq:equs-vector-form}
\begin{cases}
\mathbf{A}_t(t) = \left(\mathbf{P} - \mathbf{Q}\right) \mathbf{A}(t), & \mathbf{A}(T) = -\mathbf{1},\\
\mathbf{B}_t(t) = (\mathbf{R}-\mathbf{K}(t))\mathbf{B}(t), & \mathbf{B}(T) = -\mathbf{1},\\
\mathbf{C}_t(t) = -\mathbf{Q}\mathbf{C}(t) - \mathbf{M}(t), & \mathbf{C}(T) = \mathbf{0},\\
\mathbf{D}_t(t) = -\mathbf{Q}\mathbf{D}(t) - \mathbf{L}(t), & \mathbf{D}(T) = \mathbf{0}.
\end{cases}
\end{align}
Here, we define the following auxiliary matrices and vectors:
\begin{align*}
\mathbf{K}(t) &= \operatorname{diag}(A(t, 2)/A(t, 1), A(t, 1)/A(t, 2)) \mathbf{Q},\quad \mathbf{1} = (1, 1)^\top, \quad \mathbf{0} = (0, 0)^\top,\\
\mathbf{M}(t) &= \Big( A(t, 2)\big(B(t, 2)-B(t, 1)\big)^2,\, A(t, 1)\big(B(t, 1)-B(t, 2)\big)^2 \Big)^\top, \quad \mathbf{Q} = \begin{pmatrix} -1 & 1 \\ 1 & -1 \end{pmatrix},\\
\mathbf{L}(t) &= \frac{\gamma}{2} \bigg( \log \Big(\frac{\pi\gamma}{-\sigma_1^2 A(t, 1)}\Big),\, \log \Big(\frac{\pi \gamma}{-\sigma_2^2 A(t, 2)}\Big) \bigg)^\top,\quad \mathbf{R} = \operatorname{diag}(r_1, r_2), \quad \mathbf{P} = \operatorname{diag}(\rho_1^2-2r_1, \rho_2^2-2r_2).
\end{align*}
By standard linear SDE theory, the solutions to \eqref{eq:equs-vector-form} are determined as
\begin{align}\label{eq:equs-vector-form-solution}
\begin{cases}
\mathbf{A}(t) = -\exp \left((t-T) (\mathbf{P} - \mathbf{Q}) \right) \mathbf{1},\\
\mathbf{B}(t) = -\boldsymbol{\Phi}_B(t, T) \mathbf{1},\\
\mathbf{C}(t) = \int_{t}^{T} e^{\mathbf{Q}(s-t)} \mathbf{M}(s) \mathrm{d} s,\\
\mathbf{D}(t) = \int_{t}^{T} e^{\mathbf{Q}(s-t)} \mathbf{L}(s) \mathrm{d} s,
\end{cases}
\end{align}
Here, $\boldsymbol{\Phi}_B(t, T)$ denotes the fundamental solution matrix satisfying $\frac{\partial}{\partial t} \boldsymbol{\Phi}_B(t, T) = (\mathbf{R}-\mathbf{K}(t)) \boldsymbol{\Phi}_B(t, T)$ with $\boldsymbol{\Phi}_B(T, T) = \mathbf{I}$, where $\mathbf{I}$ denotes the $2 \times 2$ identity matrix. 
Furthermore, the matrix exponential associated with the transition intensity operator $\mathbf{Q}$ admits the explicit expression: $
e^{\mathbf{Q}(s-t)} = \frac{1}{2} \begin{pmatrix} 1 + e^{-2(s-t)} & 1 - e^{-2(s-t)} \\ 1 - e^{-2(s-t)} & 1 + e^{-2(s-t)} \end{pmatrix}$.
To facilitate the analysis, we encapsulate the model parameters into the vectors $\theta$ and $\zeta$. Identifying these with the true system parameters, we set $\theta = \zeta = (\rho_1, \rho_2, \sigma_1, \sigma_2)^\top$. The associated matrices and vector-valued functions can then be parameterized by $\theta$. Specifically, for any function $F \in \{A, B, C, D\}$, we denote $\mathbf{F}^\theta(t) = (F^\theta(t, 1), F^\theta(t, 2))^\top$. Furthermore, we define:
$\mathbf{M}^\theta(t) = \Big( A^\theta(t, 2) \big(B^\theta(t, 2) - B^\theta(t, 1)\big)^2$, $A^\theta(t, 1) \big(B^\theta(t, 1) - B^\theta(t, 2)\big)^2 \Big)^\top$,
$\mathbf{K}^\theta(t) = \operatorname{diag}(A^\theta(t, 2)/A^\theta(t, 1), A^\theta(t, 1)/A^\theta(t, 2)) \mathbf{Q}$,
$\mathbf{L}^\theta(t) = \frac{\gamma}{2} \bigg( \log \Big(\frac{\pi \gamma}{-\theta_3^2 A^\theta(t, 1)}\Big)$, $\log \Big(\frac{\pi \gamma}{-\theta_4^2 A^\theta(t, 2)}\Big) \bigg)^\top$, $\mathbf{P}^\theta = \operatorname{diag}(\theta_1^2 - 2r_1, \theta_2^2 - 2r_2)$.
Using this parameterization, the solution \eqref{eq:equs-vector-form-solution} can be rewritten in the following form
\begin{align}\label{eq:equs-vector-form-parameterize}
\left\{\begin{array}{l}
\mathbf{A}^\theta(t) = -\exp [(t-T)(\mathbf{P}^\theta-\mathbf{Q})] \mathbf{1},\\
\mathbf{B}^\theta(t) = -\mathbf{\Phi}_B^\theta(t, T) \mathbf{1},\\
\mathbf{C}^\theta(t)=\int_{t}^T e^{\mathbf{Q}(s-t)} \mathbf{M}^\theta(s) \mathrm{d} s,\\
\mathbf{D}^\theta(t) =\int_{t}^T e^{\mathbf{Q}(s-t)} \mathbf{L}^\theta(s) \mathrm{d} s.
\end{array}\right.
\end{align}
Here, $\boldsymbol{\Phi}_B^\theta(t, T)$ denotes the fundamental solution matrix satisfying $\frac{\partial}{\partial t} \boldsymbol{\Phi}_B^\theta(t, T) = (\mathbf{R} - \mathbf{K}^\theta(t)) \boldsymbol{\Phi}_B^\theta(t, T)$ with the terminal condition $\boldsymbol{\Phi}_B^\theta(T, T) = \mathbf{I}$. 
Let $\odot$ and $\oslash$ denote the Hadamard (element-wise) product and division between matrices or vectors of the same dimension, respectively. Building upon the parameterized system \eqref{eq:equs-vector-form-parameterize}, we can parameterize the value function, the q-function, and the policy in their exact forms by
\begin{align}
\mathbf{J^\theta}(t, x) =& \mathbf{A}^\theta(t) \odot \left(x \mathbf{1} + w \mathbf{B}^\theta(t)\right)^{\odot 2} + w^2 \mathbf{C}^\theta(t) + \mathbf{D}^\theta(t) + (w-z)^2 \mathbf{1},\label{eq:equs-vector-value}\\
\mathbf{q}^\zeta(t, x)
=& \left(\mathbf{A}_t^\zeta(t) + \mathbf{A}^\zeta(t) \odot \mathbf{H}_A^\zeta(t)\right) \odot \mathbf{X}^{\odot 2} + 2w \mathbf{A}^\zeta(t) \odot \left(\mathbf{B}_t^\zeta(t) - \mathbf{R}\mathbf{B}^\zeta(t)\right) \odot \mathbf{X}+ w^2 \mathbf{C}_t^\zeta(t)+ \mathbf{D}_t^\zeta(t)\notag\\
&  + \left(\mathbf{A}^\zeta(t) \odot \zeta_\sigma^{\odot 2}\right) \odot\left(a \mathbf{1}-\mathbf{m}^\zeta\right)^{\odot 2} + \mathbf{Q} \left(\mathbf{A}^\zeta(t) \odot \mathbf{X}^{\odot 2} + w^2 \mathbf{C}^\zeta(t) + \mathbf{D}^\zeta(t) + (w-z)^2 \mathbf{1}\right),\label{eq:equs-vector-q}\\
\boldsymbol{\widehat{\Pi}}_1^\zeta(a \mid t, x) =& \mathcal{N}\left(a \mid - \big(\boldsymbol{\zeta}_{\rho} \odot \mathbf{X}\big) \oslash \boldsymbol{\zeta}_{\sigma}, \; -\frac{\gamma}{2} \mathbf{1} \oslash \big[\boldsymbol{\zeta}_{\sigma}^{\odot 2} \odot \mathbf{A}^\zeta(t)\big] \right),\label{eq:equs-vector-pi}
\end{align}
where $\mathbf{H}_A^\zeta(t) = (2r_1 - \zeta_1^2, 2r_2 - \zeta_2^2)^\top$, $\mathbf{X} = x \mathbf{1} + w \mathbf{B}^\zeta(t)$, $\mathbf{\zeta}_{\rho} = (\zeta_1, \zeta_2)^\top$, and $\mathbf{\zeta}_{\sigma} = (\zeta_3, \zeta_4)^\top$, and $\mathbf{m}^\zeta=-\left(\boldsymbol{\zeta}_\rho \odot \mathbf{X}\right) \oslash \boldsymbol{\zeta}_\sigma$.
We can thereby compute the  Lagrange multiplier $\mathbf{w}$ via \eqref{eq:MV-optimal-w}
\begin{align*}
\mathbf{w}(t_0, x_0) = \big( z \mathbf{1} - x_0 \mathbf{A}(t_0) \odot \mathbf{B}(t_0) \big) \oslash \big( \mathbf{A}(t_0) \odot \mathbf{B}(t_0)^{\odot 2} + \mathbf{C}(t_0) + \mathbf{1} \big),
\end{align*}
where $\mathbf{A}(t_0) = (A(t_0, 1), A(t_0, 2))^\top$, $\mathbf{B}(t_0) = (B(t_0, 1), B(t_0, 2))^\top$, and $\mathbf{C}(t_0) = (C(t_0, 1), C(t_0, 2))^\top$ are the solutions of \eqref{eq:equs-vector-form} when $t=t_0$. The vector $\mathbf{w}(t_0, x_0)$ collects the regime-specific Lagrange multipliers. 

\begin{figure}[htb]
    \centering
\includegraphics[width=1.0\textwidth, height=0.28\textheight]{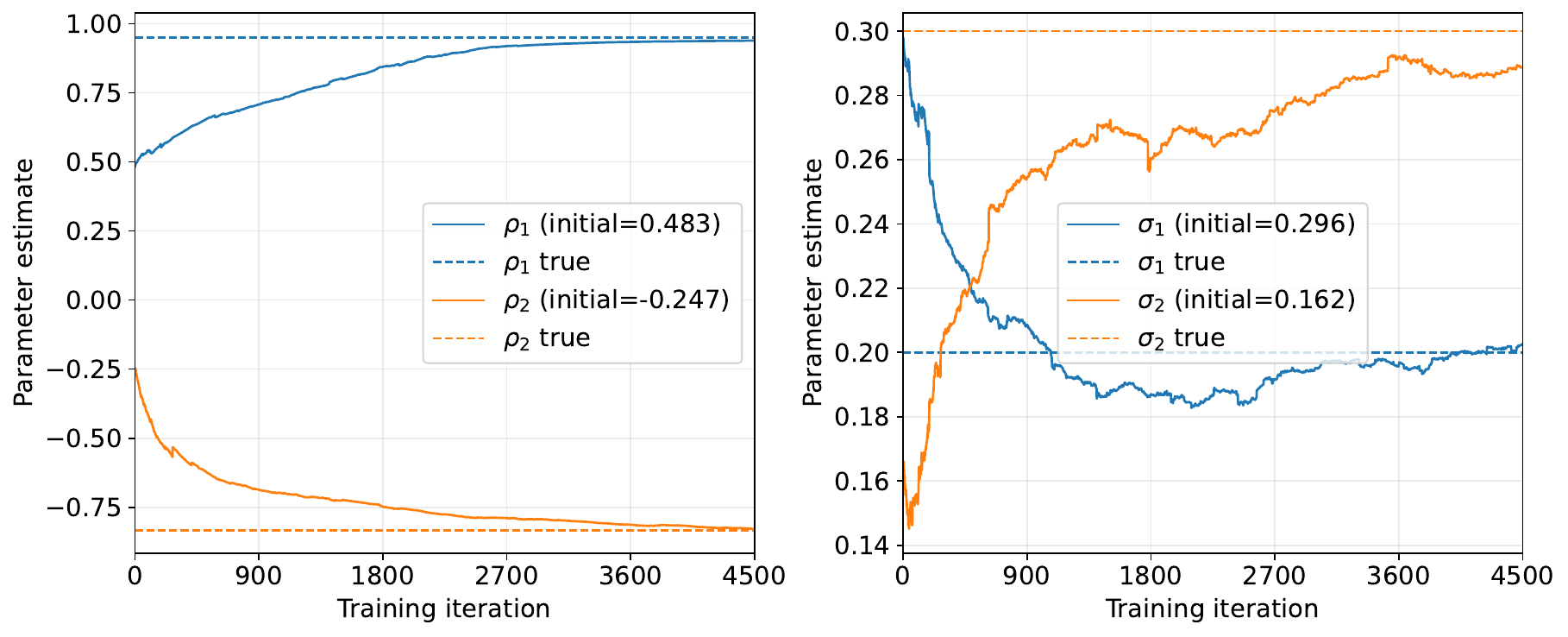}
    \caption{The panels show the convergence of parameter iterations for ($\rho_1,\rho_2,\sigma_1,\sigma_2$).}
    \label{fig:q-learning-paths-p1}
\end{figure}

To validate the theoretical results, we implement Algorithm \ref{Alg:Tsallis-q-Learning} specialized to the regime-switching MV portfolio problem, utilizing the parameterized structures in \eqref{eq:equs-vector-value}--\eqref{eq:equs-vector-pi} and the ODE system \eqref{eq:equs-vector-form-parameterize}. The investment horizon $T = 1$ (year) is discretized into $K = 25$ rebalancing epochs ($\Delta t=0.04$). For numerical integration, each decision interval $\Delta t$ is further subdivided into 20 local substeps ($h = 0.002$), over which the selected action remains fixed. We set the initial wealth $x_0 = 1$, terminal target $z = 1.4$, exploration weight $\gamma = 0.5$. For each simulated trajectory, the initial market regime is drawn uniformly from $\mathcal{M}=\{1,2\}$.
The regime-dependent market coefficients $(\mu_i, \sigma_i, r_i)$ are set to $(0.2, 0.2, 0.01)$ for the “bull” regime ($i=1$) and $(-0.2, 0.3, 0.05)$ for the “bear” regime ($i=2$). Given the risk premium $\rho_i = (\mu_i - r_i)/\sigma_i$, this yields the ground-truth vector $\vartheta_{\mathrm{true}} = (\rho_{1}, \rho_{2}, \sigma_{1}, \sigma_{2})^\top = (0.95, -0.833, 0.2, 0.3)^\top$. To eliminate notational ambiguity between the value and q-functions, we maintain a single shared parameter vector $\widehat{\vartheta}_n = (\widehat\rho_{1,n}, \widehat\rho_{2,n}, \widehat\sigma_{1,n}, \widehat\sigma_{2,n})^\top$ to concurrently parameterize the value function, q-function, and exploratory policy. The initial estimate $\widehat{\vartheta}_0$ is uniformly sampled with $\widehat\rho_{1,0} \in [0.2, 0.5]$, $\widehat\rho_{2,0} \in [-0.4, -0.1]$, and $\widehat\sigma_{1,0}, \widehat\sigma_{2,0} \in [0.15, 0.3]$. The training spans $M=4500$ iterations using Monte Carlo batches of $N_{\mathrm{paths}}=50$ independent paths. The parameter update scheme employs a root mean square-normalized per-parameter polynomial decay learning rate. The numerical results are presented in Figure \ref{fig:q-learning-paths-p1}, which plots the convergence behavior of the MV portfolio optimization problem by the offline learning algorithm within the framework of Tsallis entropy $p=1$. After sufficient iterations, the parameter estimates approach and stabilize near their corresponding ground-truth values.

\textbf{The case of $p=2$.}
We now consider the RL problem under Tsallis entropy regularization ($p > 1$), focusing on the quadratic case ($p=2$) for clarity of exposition.
Since $V_t(t, x, i; w)+r(t, i) x V_x(t, x, i; w)$ is independent of the action $a$, it can be absorbed into the normalization multiplier. For notational convenience, we continue to denote the resulting state-dependent normalizing function by $\psi(t, x, i)$, while retaining the regime-switching term $\sum_{j=1}^L q_{i j} V(t, x, j; w)$ explicitly. Besides, we define $m(t,x,i) = -\frac{\rho(t,i)}{\sigma(t,i)}[x+wB(t,i)]$, $\lambda(t,i) = -A(t,i)\sigma^2(t,i)>0$. For $p=2$, the optimal policy in \eqref{eq:MV-optimal-policy-m1m2m3-p2} becomes
\begin{align}\label{eq:MV-policy-p2-original}
\widehat{\pi}_2(a\mid t,x,i) =& \frac{1}{2\gamma} \Bigg[ 2A(t,i)\rho(t,i)\sigma(t,i)[x+wB(t,i)]a +A(t,i)\sigma^2(t,i)a^2 +\sum_{j=1}^{L}q_{ij}V(t,x,j;w) +\psi(t,x,i) \Bigg]_+ \nonumber\\
=& \frac{1}{2\gamma} \left[ \Gamma(t,x,i) -\lambda(t,i)[a-m(t,x,i)]^2 \right]_+,
\end{align}
where $\psi(t,x,i)$ is determined by $\int_{\mathbb{R}} \widehat{\pi}_2(a \mid t, x, i) \mathrm{d}a = 1$, and $\Gamma(t,x,i) = \sum_{j=1}^{L}q_{ij}V(t,x,j;w) +\psi(t,x,i) -A(t,i)\rho^2(t,i)[x+wB(t,i)]^2$.
Since $\widehat{\pi}_2(\cdot\mid t,x,i)$ is a probability density, the normalization condition requires $\Gamma(t,x,i)>0$. Define the radius $R(t, x, i)=\big(\frac{\Gamma(t, x, i)}{\lambda(t, i)}\big)^{1/2}>0$.
Then we have
\begin{align}\label{eq:MV-policy-p2-centered}
\widehat{\pi}_2(a\mid t,x,i) = \frac{\lambda(t,i)}{2\gamma} \left[ R^2(t, x, i)-[a-m(t,x,i)]^2 \right]_+.
\end{align}
In particular, its support is given by $
\operatorname{supp} \widehat{\pi}_2(\cdot\mid t,x,i) = [m(t,x,i)-R(t, x, i),\,m(t,x,i)+R(t, x, i)]$.
Thus, although the admissible action space is $\mathcal A=\mathbb R$, the optimal policy under the quadratic Tsallis regularization is endogenously compactly supported. We next determine $R(t, x, i)$ from the normalization condition. By \eqref{eq:MV-policy-p2-centered},
$1 = \int_{\mathbb R} \widehat{\pi}_2(a\mid t,x,i)\mathrm{d} a = \frac{\lambda(t,i)}{2\gamma} \int_{m(t,x,i)-R(t, x, i)}^{m(t,x,i)+R(t, x, i)} \left[ R^2(t, x, i)-[a-m(t,x,i)]^2 \right]\mathrm{d} a. $
Introducing the change of variable $u=a-m(t,x,i)$, we obtain
\begin{align*}
1 = \frac{\lambda(t,i)}{2\gamma} \int_{-R(t, x, i)}^{R(t, x, i)} \left[R^2(t, x, i)-u^2\right]\mathrm{d}u = \frac{\lambda(t,i)}{2\gamma} \left[ 2R^3(t,x,i)-\frac{2}{3}R^3(t,x,i) \right] = \frac{2\lambda(t,i)}{3\gamma}R^3(t,x,i).
\end{align*}
Therefore,
\begin{align}\label{eq:MV-p2-radius}
R(t, x, i) = \left( \frac{3\gamma}{2\lambda(t,i)} \right)^{1/3} = \left( \frac{3\gamma} {-2A(t,i)\sigma^2(t,i)} \right)^{1/3}.
\end{align}
Hence, the radius is independent of the state variable $x$, and we shall henceforth write it as $R(t,i)$.
Substituting \eqref{eq:MV-p2-radius} into \eqref{eq:MV-policy-p2-centered} yields the explicit optimal policy
\begin{align}\label{eq:MV-optimal-policy-p2-explicit}
\widehat{\pi}_2(a\mid t,x,i) = \frac{-A(t,i)\sigma^2(t,i)}{2\gamma} \left[ R^2(t,i) - \left( a+ \frac{\rho(t,i)}{\sigma(t,i)} [x+wB(t,i)] \right)^2 \right]_+,
\end{align}
where $R(t,i)$ is given by \eqref{eq:MV-p2-radius}.
It is worth noting that retaining the regime-switching term $\sum_{j=1}^{L}q_{ij}V(t,x,j;w)$ in \eqref{eq:MV-policy-p2-original} does not alter the center or the radius after the normalization constraint is imposed. Indeed, combining \eqref{eq:MV-policy-p2-original} and \eqref{eq:MV-optimal-policy-p2-explicit} yields
\begin{align}\label{eq:MV-p2-psi-explicit}
\psi(t,x,i) = &\lambda(t,i)R^2(t,i) -\sum_{j=1}^{L}q_{ij}V(t,x,j;w) +A(t,i)\rho^2(t,i)[x+wB(t,i)]^2\nonumber\\
= &\left(\frac{3\gamma}{2}\right)^{2/3} [-A(t,i)\sigma^2(t,i)]^{1/3} -\sum_{j=1}^{L}q_{ij}V(t,x,j;w)+A(t,i)\rho^2(t,i)[x+wB(t,i)]^2.
\end{align}
Hence, the regime-switching contribution is absorbed through the normalizing function after imposing the probability constraint, while the shape of the optimal density is determined by its center $m(t,x,i)$ and radius $R(t,i)$. We now derive the first two moments of $\widehat{\pi}_2(\cdot\mid t,x,i)$. Since the density in \eqref{eq:MV-policy-p2-centered} is symmetric around $m(t,x,i)$, we have $\int_{\mathbb R} [a-m(t,x,i)] \widehat{\pi}_2(a\mid t,x,i)\mathrm{d} a =0$.
Therefore, the first moment is
\begin{align}\label{eq:MV-first-moment-p2}
\mathbb E_{\widehat{\pi}_2}[a] = m(t,x,i) = -\frac{\rho(t,i)}{\sigma(t,i)} [x+wB(t,i)].
\end{align}
To obtain the second moment, we first compute the second centered moment. Using again the change of variable $u=a-m(t,x,i)$, we have
$\mathbb E_{\widehat{\pi}_2} \left[ (a-m(t,x,i))^2 \right] = \frac{\lambda(t,i)}{2\gamma} \int_{-R(t,i)}^{R(t,i)} u^2 \left[ R^2(t,i)-u^2 \right]\mathrm{d} u  = \frac{2\lambda(t,i)}{15\gamma}R^5(t,i). $
From the normalization identity $\lambda(t,i)R^3(t,i) = \frac{3\gamma}{2}$,
it follows that $
\operatorname{Var}_{\widehat{\pi}_2}[a] = \mathbb E_{\widehat{\pi}_2} \left[ (a-m(t,x,i))^2 \right] = \frac{R^2(t,i)}{5}$.
Consequently,
\begin{align}\label{eq:MV-second-moment-p2-explicit}
\mathbb E_{\widehat{\pi}_2}[a^2] &= m^2(t,x,i) + \frac{R^2(t,i)}{5}= \frac{\rho^2(t,i)}{\sigma^2(t,i)} [x+wB(t,i)]^2 + \frac{1}{5} \left( \frac{3\gamma} {-2A(t,i)\sigma^2(t,i)} \right)^{2/3}.
\end{align}
For subsequent use in the Tsallis entropy term, we also record the squared $L^2$-norm of the optimal density. From \eqref{eq:MV-policy-p2-centered}, we have
\begin{align}\label{eq:MV-p2-L2}
\int_{\mathbb R} \widehat{\pi}_2^2(a\mid t,x,i)\mathrm{d} a
&= \frac{\lambda^2(t,i)}{4\gamma^2} \int_{-R(t,i)}^{R(t,i)} \left[ R^2(t,i)-u^2 \right]^2\mathrm{d} u = \frac{\lambda^2(t,i)}{4\gamma^2} \frac{16}{15}R^5(t,i)= \frac{3}{5R(t,i)}.
\end{align}
Equations \eqref{eq:MV-first-moment-p2}, \eqref{eq:MV-second-moment-p2-explicit}, and \eqref{eq:MV-p2-L2} provide explicit expressions for all policy-dependent quantities entering the HJB equation. In particular, the regime-switching term $\sum_{j=1}^{L}q_{ij}V(t,x,j;w)$ affects the normalizing function $\psi(t,x,i)$ through \eqref{eq:MV-p2-psi-explicit}, but does not alter the center, support radius, or the first two moments of the normalized policy. Moreover, the dependence of the first two moments on the wealth state $x$ occurs exclusively through the shifted state $x+w B(t, i)$, while the exploration component $R^2(t, i)/5$ is independent of $x$. This structure permits coefficient-wise matching in powers of $x+wB(t,i)$ in the HJB equation and thereby leads to a closed system of ordinary differential equations for the coefficient functions $A$, $B$, $C$, and $D$.
Substituting the first and second moments $\Ex_{\widehat{\pi}_2}[a]$ and $\Ex_{\widehat{\pi}_2}[a^2]$ of the optimal policy $\widehat{\pi}_2(a \mid t, x, i)$ into the HJB equation \eqref{eq:MV-HJB} and matching the coefficients of
$[x+wB(t,i)]^2$, $[x+wB(t,i)]$, and the remaining state-independent terms yield
\begin{align}\label{eq:MV-ODEs-p2}
\begin{cases}
A_t(t,i)=\big[\rho^2(t,i)-2r(t,i)\big]A(t,i)-\sum_{j=1}^2q_{ij}A(t,j), &A(T,i)=-1,\\
B_t(t,i)=r(t,i)B(t,i)-\frac{1}{A(t,i)}\sum_{j=1}^2q_{ij}A(t,j)\big[B(t,j)-B(t,i)\big],&B(T,i)=-1,\\
C_t(t,i)=-\sum_{j=1}^2q_{ij}C(t,j)-\sum_{j=1}^2q_{ij}A(t,j)\big[B(t,j)-B(t,i)\big]^2, &C(T,i)=0,\\
D_t(t,i)=-\sum_{j=1}^2q_{ij}D(t,j)-\gamma+\frac{3}{5}\lambda(t,i)R^2(t,i),
&D(T,i)=0.
\end{cases}
\end{align}

We next characterize the corresponding q-function defined in
\eqref{eq:MV-q-func}. Substituting the quadratic value function and using
\eqref{eq:MV-ODEs-p2}, we obtain
\begin{align}\label{eq:MV-q-function-p2}
q^{\widehat{\pi}_2}(t,x,i,a;w)
=&J^{\widehat{\pi}_2}_t(t,x,i;w)+\big[r(t,i)x+\rho(t,i)\sigma(t,i)a\big]J^{\widehat{\pi}_2}_x(t,x,i;w)
\nonumber\\
&+\frac{1}{2}\sigma^2(t,i)a^2J^{\widehat{\pi}_2}_{xx}(t,x,i;w)+\sum_{j=1}^2q_{ij}J^{\widehat{\pi}_2}(t, x, j; w)\nonumber\\
=&-\lambda(t,i)\big[a-m(t,x,i)\big]^2
-\gamma+\frac{3}{5}\lambda(t,i)R^2(t,i).
\end{align}
For subsequent analysis, we define the policy-averaged q-function as $\overline q^{\widehat{\pi}_2}(t,x,i;w) := \int_{\mathbb R}q^{\widehat{\pi}_2}(t,x,i,a;w)\widehat{\pi}_2(a\mid t,x,i)\mathrm{d} a$. By evaluating this expectation over the action-dependent q-function \eqref{eq:MV-q-function-p2} and substituting the variance $\Ex_{\widehat{\pi}_2}\big[(a-m(t,x,i))^2\big]=\frac{R^2(t,i)}{5}$, we establish the following identity strictly for the policy-averaged q-function:
\begin{equation}\label{eq:MV-qbar-p2}
    \overline q^{\,\widehat{\pi}_2}(t,x,i;w) = -\gamma+\frac{2}{5}\lambda(t,i)R^2(t,i) = -\gamma+\gamma\int_{\mathbb R}\widehat{\pi}_2^2(a\mid t,x,i)\mathrm{d} a.
\end{equation}

For a compact vector representation, define, for any
$F\in\{A,B,C,D\}$, $\mathbf F(t)=(F(t,1),F(t,2))^\top$, and let $\boldsymbol{\rho}(t)=(\rho(t,1),\rho(t,2))^\top$, $\boldsymbol{\sigma}(t)=(\sigma(t,1),\sigma(t,2))^\top$.
Then \eqref{eq:MV-ODEs-p2} can be written as
\begin{align}\label{eq:MV-ODEs-vector-p2}
\begin{cases}
\mathbf A_t(t)=\big[\mathbf P(t)-Q\big]\mathbf A(t),&\mathbf A(T)=-\mathbf{1},\\
\mathbf B_t(t)=\mathbf R_f(t)\mathbf B(t)-\mathbf N_B(t),&\mathbf B(T)=-\mathbf{1},\\
\mathbf C_t(t)=-Q\mathbf C(t)-\mathbf M(t),&\mathbf C(T)=\mathbf 0,\\
\mathbf D_t(t)=-Q\mathbf D(t)-\mathbf L_D(t),&\mathbf D(T)=\mathbf 0,\end{cases}
\end{align}
where $\mathbf P(t)=\operatorname{diag}\big(\rho^2(t,1)-2r(t,1), \rho^2(t,2)-2r(t,2)
\big)$, $\mathbf R_f(t)=\operatorname{diag}\big(r(t,1),r(t,2)\big)$,
 $\boldsymbol{\lambda}(t)=-\mathbf A(t)\odot\boldsymbol{\sigma}(t)^{\odot 2}$, $\mathbf N_B(t)=\left[Q\big(\mathbf A(t)\odot\mathbf B(t)\big)-\mathbf B(t)\odot Q\mathbf A(t)\right]\oslash \mathbf A(t)$,
$\mathbf M(t)=Q\big(\mathbf A(t)\odot\mathbf B(t)^{\odot 2}\big)-2\mathbf B(t)\odot Q\big(\mathbf A(t)\odot \mathbf B(t)\big)+\mathbf B(t)^{\odot 2}\odot Q\mathbf A(t)$,
 $\mathbf L_D(t)=\gamma\mathbf{1}-\frac{3}{5}\boldsymbol{\lambda}(t)\odot\mathbf R_\pi(t)^{\odot 2}$, $\mathbf R_\pi(t)=\left[\frac{3\gamma}{2}\mathbf{1}\oslash\boldsymbol{\lambda}(t)\right]^{\odot 1/3}$. Since $A(t,i)<0$ and $\sigma(t,i)>0$, $\boldsymbol{\lambda}(t)$ and $\mathbf R_\pi(t)$ are component-wise positive.
By standard linear ODE theory, the coefficient system admits the following representation
\begin{align}\label{eq:MV-ODEs-solution-p2}
\begin{cases}
\mathbf A(t)=-\boldsymbol{\Phi}_A(t,T)\mathbf{1},\\
\mathbf B(t)=-\boldsymbol{\Phi}_B(t,T)\mathbf{1},\\
\mathbf C(t)=
\int_t^T e^{Q(s-t)}\mathbf M(s)\mathrm{d} s,\\
\mathbf D(t)=\int_t^T e^{Q(s-t)}\mathbf L_D(s)\mathrm{d} s.
\end{cases}
\end{align}
Here, $\boldsymbol{\Phi}_A(t,T)$ denotes the fundamental solution matrix satisfying $\frac{\partial}{\partial t}\boldsymbol{\Phi}_A(t,T)
=[\mathbf P(t)-Q]\boldsymbol{\Phi}_A(t,T)$ with  $\boldsymbol{\Phi}_A(T,T)=\mathbf I$ and $\boldsymbol{\Phi}_B(t,T)$ denotes the fundamental solution matrix satisfying
$\frac{\partial}{\partial t}\boldsymbol{\Phi}_B(t,T)=[\mathbf R_f(t)-\mathbf K_A(t)]\boldsymbol{\Phi}_B(t,T)$ with $\boldsymbol{\Phi}_B(T,T)=\mathbf I$,
where $\mathbf K_A(t)=\operatorname{diag}(\mathbf A(t))^{-1}
Q\operatorname{diag}(\mathbf A(t))-\operatorname{diag}\left(Q\mathbf A(t)\oslash\mathbf A(t)\right)$, and $\mathbf N_B(t)=\mathbf K_A(t)\mathbf B(t)$.

We next express the optimal policy and its moments in vector form.
Set $\mathbf X(t,x)=x\mathbf{1}+w\mathbf B(t)$, $\mathbf m(t,x)=-\big[\boldsymbol{\rho}(t)\odot\mathbf X(t,x)\big]\oslash
\boldsymbol{\sigma}(t)$.
Then
\begin{align}\label{eq:MV-vector-policy-p2}
\boldsymbol{\widehat{\Pi}}_2(a\mid t,x)=\frac{1}{2\gamma}\boldsymbol{\lambda}(t)
\odot\left[\mathbf R_\pi(t)^{\odot 2}-\big(a\mathbf{1}-\mathbf m(t,x)\big)^{\odot 2}\right]_+.
\end{align}
Accordingly, $\mathbf E_{\boldsymbol{\widehat{\Pi}}_2}[a]=\mathbf m(t,x)$,
$\mathbf E_{\boldsymbol{\widehat{\Pi}}_2}[a^2]=\mathbf m(t,x)^{\odot 2}+\frac{1}{5}\mathbf R_\pi(t)^{\odot 2}$, $\int_{\mathbb R}\boldsymbol{\widehat{\Pi}}_2(a\mid t,x)^{\odot 2}\mathrm{d} a
=\frac{3}{5}\mathbf{1}\oslash\mathbf R_\pi(t)$.
Here, the positive-part operator and all fractional powers are understood componentwise.

In the numerical implementation below, we further specialize to time-homogeneous regime-dependent market coefficients (i.e., $\rho(t, i) \equiv \rho_i$, $r(t, i) \equiv r_i$, and $\sigma(t, i) \equiv \sigma_i$). 
We next introduce separate parameter copies for the value
function, the q-function, and the policy function. Let $\vartheta^\star=(\rho_1,\rho_2,\sigma_1,\sigma_2)^\top$ denote the true model parameter. We use
$\theta$, $\zeta$, and $\chi$ as functionally distinct parameter
copies for the value
function, the q-function, and the policy function parametrizations,
respectively. They are not three separate sets of economic parameters; instead, they are three functionally distinct copies of the same structural market parameter vector. Although the three parameter copies share the same population target $\theta^\star=\zeta^\star=\chi^\star=\vartheta^\star$,
their current estimates $\theta$, $\zeta$, $\chi$ are treated as functionally distinct and are updated separately during learning.
For each parameter copy
$\eta\in\{\theta,\zeta,\chi\}$, let
$\boldsymbol{\eta}_{\rho}=(\eta_1,\eta_2)^\top$, $\boldsymbol{\eta}_{\sigma}
=(\eta_3,\eta_4)^\top$.
We denote by $(\mathbf A^\eta,\mathbf B^\eta,\mathbf C^\eta,\mathbf D^\eta)$
the solution of the coefficient system
\eqref{eq:MV-ODEs-vector-p2} corresponding
to $\eta$. Define $\mathbf X^\eta(t,x)=x\mathbf{1}+w\mathbf B^\eta(t)$,
$\boldsymbol{\lambda}^\eta(t)=-\mathbf A^\eta(t)\odot(\boldsymbol{\eta}_\sigma)^{\odot 2}$,
$\mathbf m^\eta(t,x)=-\big[\boldsymbol{\eta}_\rho\odot\mathbf X^\eta(t,x)\big]\oslash\boldsymbol{\eta}_\sigma$,
$\mathbf R_\pi^\eta(t)=\left[\frac{3\gamma}{2}\mathbf{1}\oslash\boldsymbol{\lambda}^\eta(t)\right]^{\odot 1/3}$, $\widehat{\Pi}_2^\eta(a \mid t, x)=\frac{1}{2 \gamma} \boldsymbol{\lambda}^\eta(t) \odot\left[\left(\mathbf{R}_\pi^\eta(t)\right)^{\odot 2}-\left(a \mathbf{1}-\mathbf{m}^\eta(t, x)\right)^{\odot 2}\right]_{+}$.
The parameterized value function is therefore
\begin{align}\label{eq:MV-V-parameter-p2}
\mathbf V^\eta(t,x;w)=\mathbf A^\eta(t)\odot(\mathbf X^\eta(t,x))^{\odot2}+w^2\mathbf C^\eta(t)+\mathbf D^\eta(t)
+(w-z)^2\mathbf1, \,\,\eta \in\{\theta, \zeta, \chi\}.
\end{align}
For the q-function parameter $\zeta$, the corresponding action-dependent q-function is
\begin{align}\label{eq:MV-q-parameter-p2}
\mathbf q^\zeta(t,x,a;w)
=-\boldsymbol{\lambda}^\zeta(t)\odot\big[a\mathbf{1}-\mathbf m^\zeta(t,x)
\big]^{\odot 2}-\gamma\mathbf{1}+\frac{3}{5}\boldsymbol{\lambda}^\zeta(t)\odot\big(\mathbf R_\pi^\zeta(t)\big)^{\odot 2}.
\end{align}
The critic value-function approximation is \(\mathbf V^\theta\), whereas \(\mathbf V^\zeta\) and \(\mathbf V^\chi\) are auxiliary structural representations used in constructing the q-function and policy, respectively.
The q-function parameterization is restricted to the exact structural family induced by the coefficient system corresponding to $\zeta$. The numerical experiment uses a structure-preserving parametric implementation rather than an unrestricted function approximator. For later use, we define the self-averaged q-function associated with the structural parameter $\zeta$ by
\begin{align*}
\overline{\mathbf q}^{\,\zeta}(t,x;w)=-\gamma\mathbf{1}+\gamma\int_{\mathbb R}\boldsymbol{\widehat{\Pi}}_2^\zeta(a\mid t,x)^{\odot 2}\mathrm{d}a=-\gamma\mathbf{1}+\frac{2}{5}\boldsymbol{\lambda}^\zeta(t)\odot\big(\mathbf {R}_\pi^\zeta(t)\big)^{\odot 2}.
\end{align*}
Thus, the HJB identity provides an exact expression for the self-averaged q-function obtained by averaging $\mathbf{q}^\zeta$ under the policy $\widehat{\Pi}_2^\zeta$ induced by the same structural parameter $\zeta$. During learning, however, the actor parameter $\chi$ need not coincide with $\zeta$; hence, an average of $\mathbf{q}^\zeta$ under the current actor policy $\widehat{\Pi}_2^\chi$ should not be identified with $\overline{\mathbf{q}}^\zeta$ unless $\chi=\zeta$.

As in the preceding derivation, the action-independent term $V_t^\chi+r(t, i) x V_x^\chi$ is absorbed into the normalizing vector, while the regime-switching contribution $Q \mathbf{V}^\chi$ is retained explicitly. This separation is purely notational: both terms are independent of the action and could in principle be absorbed into the normalization multiplier. We retain $Q \mathbf{V}^\chi$ only to preserve the regime-switching structure of the original Hamiltonian representation. Finally, to explicitly retain the regime-switching contribution in the original policy representation for parameter $\chi$, we define 
$
\mathbf K_1^\chi(t,x)=2\mathbf A^\chi(t)\odot\mathbf X^\chi(t,x)\odot\boldsymbol{\chi}_\rho\odot\boldsymbol{\chi}_\sigma$, $
\mathbf K_2^\chi(t)=\mathbf A^\chi(t)\odot(\boldsymbol{\chi}_\sigma)^{\odot 2}$,
and $\boldsymbol{\psi}^\chi(t,x)=\boldsymbol{\lambda}^\chi(t)\odot\big(\mathbf R_\pi^\chi(t)\big)^{\odot 2}-Q\mathbf V^\chi(t,x;w)+\mathbf A^\chi(t)\odot(\boldsymbol{\chi}_\rho)^{\odot 2}\odot\big(\mathbf X^\chi(t,x)\big)^{\odot 2}$.
The policy can then be written in the original Hamiltonian form as
\begin{align}\label{eq:MV-policy-parameter-p2}
\boldsymbol{\widehat{\Pi}}_2^\chi(a\mid t,x)
&=\frac{1}{2\gamma}\left[\mathbf K_1^\chi(t,x)a+\mathbf K_2^\chi(t)a^2+Q\mathbf V^\chi(t,x;w)
+\boldsymbol{\psi}^\chi(t,x)\right]_+ \nonumber\\
&=\frac{1}{2\gamma}
\boldsymbol{\lambda}^\chi(t)\odot\left[\big(\mathbf R_\pi^\chi(t)\big)^{\odot 2}-\big(a\mathbf{1}-\mathbf m^\chi(t,x)\big)^{\odot 2}\right]_+.
\end{align}
Although $Q \mathbf{V}^\chi$ appears explicitly in the pre-normalized Hamiltonian representation, it is an action-independent term and is exactly offset by the corresponding component of the normalizing vector $\boldsymbol{\psi}^\chi$. Hence, conditional on the coefficient functions $\mathbf{A}^\chi$ and $\mathbf{B}^\chi$, the explicit additive term $Q \mathbf{V}^\chi$ does not contribute an additional explicit term to the center, support radius, or moments of the normalized density. Regime switching nevertheless affects these quantities indirectly through the coupled coefficient system determining $\mathbf{A}^\chi$ and $\mathbf{B}^\chi$.

\begin{figure}[htb]
    \centering
    \includegraphics[width=1.00\textwidth, height=0.28\textheight]{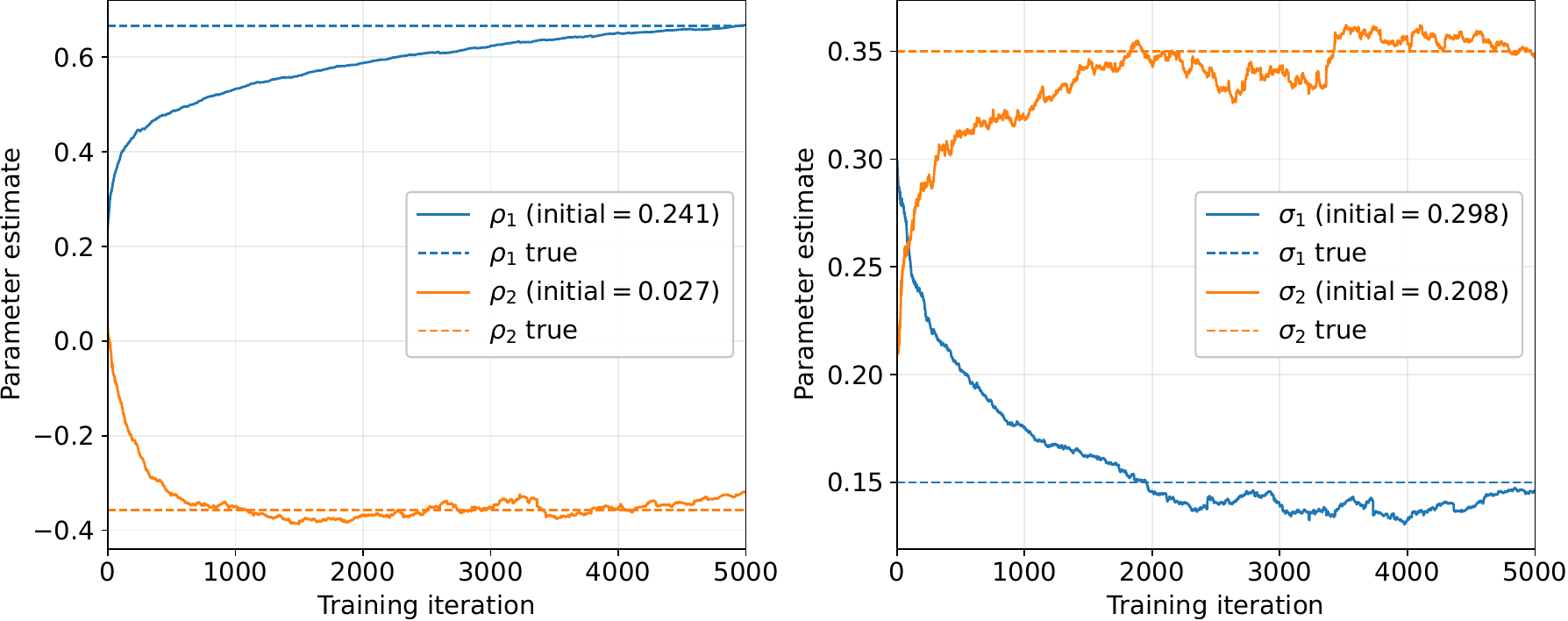}
    \caption{The panels show the convergence of parameter iterations for ($\rho_1,\rho_2,\sigma_1,\sigma_2$).}
    \label{fig:q-learning-paths-p2}
\end{figure}

To validate the theoretical framework for the Tsallis entropy case ($p=2$), we implement Algorithm \ref{Alg:Tsallis-q-Learning} specialized to the regime-switching continuous-time MV portfolio problem. The investment horizon $T = 1$ (year) is discretized into $K = 25$ rebalancing epochs ($\Delta t=0.04$). For numerical integration, each decision interval $\Delta t$ is further subdivided into 20 local substeps ($h = 0.002$), over which the selected action remains fixed. We set the initial wealth $x_0 = 1$, terminal target $z = 1.4$, and the temperature parameter $\gamma = 0.5$. The Markov chain generator $Q$ is configured with transition rates $q_{12} = 1.8$ and $q_{21} = 2.0$
. For each simulated trajectory, the initial market regime is drawn uniformly from $\mathcal{M}=\{1, 2\}$.
The regime-dependent market coefficients $(\mu_i, \sigma_i, r_i)$ are set to $(0.12, 0.15, 0.02)$ for the “bull” regime ($i=1$) and $(-0.10, 0.35, 0.025)$ for the “bear” regime ($i=2$). Following the analytical definition of the risk premium $\rho_i = (\mu_i - r_i)/\sigma_i$, this configuration mathematically rigorously yields the ground-truth parameter vector $\vartheta_{\mathrm{true}} = (\rho_{1}, \rho_{2}, \sigma_{1}, \sigma_{2})^\top = (0.6667, -0.3571, 0.15, 0.35)^\top$. To eliminate notational ambiguity, we maintain a single shared parameter vector $\widehat{\vartheta}_n = (\widehat\rho_{1,n}, \widehat\rho_{2,n}, \widehat\sigma_{1,n}, \widehat\sigma_{2,n})^\top$ to concurrently parameterize the value function, q-function, and exploratory policy. The initial estimate $\widehat{\vartheta}_0$ is uniformly sampled with $\widehat\rho_{1,0} \in [0.10, 0.25]$, $\widehat\rho_{2,0} \in [-0.05, 0.10]$, and $\widehat\sigma_{1,0}, \widehat\sigma_{2,0} \in [0.20, 0.30]$. Notably, the theoretical action space is maintained as unconstrained ($\mathcal{A} = \mathbb{R}$) without any exogenous truncation, as the optimal Tsallis policy natively generates a probability density with compact support $[m_i^\zeta - R_{\pi,i}^\zeta, m_i^\zeta + R_{\pi,i}^\zeta]$. The training spans $M=5000$ iterations using Monte Carlo batches of $N_{\mathrm{paths}}=50$ independent paths. The parameter update scheme employs a root mean square-normalized per-parameter polynomial decay learning rate. The numerical results are presented in Figure \ref{fig:q-learning-paths-p2}, which plots the convergence behavior of the MV portfolio optimization problem by the proposed q-learning algorithm within the framework of Tsallis entropy $p=2$. The numerical trajectories exhibit clear convergence toward the ground-truth parameter values.

\section{Conclusion}
Building upon prior research on regime-switching models and recent advances in RL for stochastic control, this paper extends continuous-time $q$-learning theory by incorporating Tsallis entropy as the regularizer.
Conventional RL predominantly utilizes Shannon entropy. The more general Tsallis entropy is often avoided because it yields optimal policy distributions outside standard probabilistic families, which complicates subsequent algorithm design. This work systematically resolves this analytical hurdle.
Furthermore, existing continuous-time regime-switching RL algorithms are largely confined to the EMV framework. This restricts their applicability to broader stochastic control paradigms, such as LQ problems. To address this limitation, we investigate continuous-time $q$-learning for Markov regime-switching systems under Tsallis entropy. Our objective is to develop a robust, universally applicable RL framework capable of modeling sudden environmental shifts and uncertainties. Specifically, we establish a martingale characterization of the q-function with respect to Tsallis entropy for a Markov regime-switching system. Based on this characterization, we design two classes of q-learning algorithms, categorized by whether the Lagrange multiplier can be explicitly derived. When the Lagrange multiplier is intractable, we propose two distinct policy update schemes: an update rule based on \eqref{eq:update-rule-F} (Algorithm \ref{Alg:Tsallis-q-Learning-normalizing-unavailable-F}) and an update rule based on KL divergence \eqref{eq:KL-update}. 
We apply these algorithms to the continuous-time EMV portfolio optimization problem in a regime-switching market. 

\vspace{0.25in}


%









\end{document}